\documentclass[a4,12pt]{amsart}

\usepackage{mathrsfs}

\usepackage{caption}

\usepackage{amsthm}
\usepackage{hyperref}
\usepackage{amssymb}
\usepackage{latexsym}
\usepackage{amsmath}
\usepackage{amsfonts}
\usepackage{mathabx}
\usepackage[dvips]{graphicx}
\usepackage[utf8]{inputenc}

\usepackage[LGR,T1]{fontenc}%
\usepackage[french,english]{babel}
\usepackage{url}
\usepackage{enumerate}
\usepackage{hyperref}
\usepackage{color}

\usepackage{fancyhdr}
\usepackage{bbm}
\usepackage{tikz}

\newtheorem{Theorem}{Theorem}[part]
\newtheorem{Definition}{Definition}[part]
\newtheorem{Proposition}{Proposition}[part]

\newtheorem{Assumption}{Assumption}[part]
\newtheorem{Lemma}{Lemma}[part]
\newtheorem{Corollary}{Corollary}[part]
\newtheorem{Remark}{Remark}[part]

\makeatletter \@addtoreset{equation}{section}

\@addtoreset{Definition}{section}

\@addtoreset{Theorem}{section}

\@addtoreset{Proposition}{section}

\@addtoreset{Property}{section}

\@addtoreset{Assumption}{section}

\@addtoreset{Corollary}{section}

\@addtoreset{Lemma}{section}

\@addtoreset{Remark}{section}

\@addtoreset{Example}{section}

\newcommand{\cC}{\mathcal{C}}
\newcommand{\cD}{\mathcal{D}}

\newcommand{\cF}{\mathcal{F}}

\newcommand{\cK}{\mathcal{K}}

\newcommand{\cS}{\mathcal{S}}

\renewcommand{\P}{\mathbb{P}}

\newcommand{\R}{\mathbb{R}}

\newcommand{\dZ}
    {\ensuremath{\delta Z }}
    
\newcommand{\dY}
    {\ensuremath{\delta Y }}

\def \proof{{\noindent \bf Proof. }}
\def \eproof{\hbox{ }\hfill$\Box$}

\newcommand{\ud}{\mathrm{d}}

\newcommand{\1}{{\bf 1}}
    
\newcommand{\set}[1]
    {\ensuremath{\{ #1 \}}}
    
\newcommand{\HP}[1] 
    {\ensuremath{\mathscr{H}^{#1}}}

\newcommand{\SP}[1]{\ensuremath{\mathscr{S}^{#1}}}

\newcommand{\esp}[1]{\ensuremath{\mathbb{E} \!\! \left[#1\right] }}

\newcommand{\EFp}[2]
    {\ensuremath{
     \mathbb{E}_{#1}\!\!\left[#2\right] }}
     
\renewcommand{\Xi}[1]{X_{i #1}}

\newcommand{\df}{\delta f}

\newcommand{\BP}[1]{\ensuremath{\mathscr{B}^{#1}}}
\newcommand{\NB}[2]{\ensuremath{\left\| #2 \right\|_{\mathscr{B}^{#1}}} }
\newcommand{\NL}[2]{\ensuremath{\left\| #2 \right\|_{\mathscr{L}^{#1}}} }

\newcommand{\bone}{\mathbf{1}}
\newcommand{\EE}{\mathbb{E}}

\begin{document}
\title{Reflected BSDEs in non-convex domains}
\author{Jean-François CHASSAGNEUX, Sergey NADTOCHIY, Adrien RICHOU}
\address{UFR de Math\'ematiques \& LPSM, Université de Paris.}
\email{chassagneux@lpsm.paris}
\address{Department of Applied Mathematics, Illinois Institute of Technology, Chicago, IL 60616.}
\email{{snadtochiy@iit.edu}}
\address{Université de Bordeaux, IMB, UMR 5251, F-33400 Talence, France.}
\email{adrien.richou@math.u-bordeaux.fr}
\footnotetext{Authors would like to thank Marc Arnaudon for the enlightening discussions about martingales on manifolds. The authors also thank the Illinois Institute of Technology for hosting the meetings during which this research was initiated. Partial support from the NSF CAREER grant 1855309 is acknowledged.}

\begin{abstract}
This paper establishes the well-posedness of reflected backward stochastic differential equations in the non-convex domains that satisfy a weaker version of the star-shaped property. The main results are established (i) in a Markovian framework with H\"older-continuous generator and terminal condition and (ii) in a general setting under a smallness assumption on the input data. We also investigate the connections between this well-posedness result and the theory of martingales on manifolds.
\end{abstract}

\maketitle



\section{Introduction}


Backward stochastic differential equations (BSDEs), originally introduced in \cite{Bismut-76} and fully developed in \cite{Peng-91,Pardoux-Peng-92}, can be viewed as the probabilistic analogues of semi-linear partial differential equations (PDEs). In particular, BSDEs are used to describe the solutions of stochastic control problems (see, among many others, \cite{Peng-93,ElKaroui-Peng-Quenez-97,Hamadene-Lepeltier-Peng-97}). If the control variable of such an optimization problem has a discrete component -- e.g., an option to switch the state process to a different regime or to terminate the process and obtain an instantaneous payoff -- then, the associated PDE obtains a free-boundary feature and the associated BSDE becomes reflected: i.e., its solution lives inside a given domain and is reflected at the boundary of this domain. The theory of reflected BSDEs in dimension one, i.e. when the reflected process is one-dimensional, is well developed in a very high generality: see, e.g. \cite{el1997reflected,cvitanic1996backward,dumitrescu2016generalized,Grigorova-Imkeller-Offen-Ouknine-Quenez-17,Grigorova-Imkeller-Ouknine-Quenez-18}.
However, the multidimensional case presents significant additional challenges (e.g., due to the lack of the comparison principle), and, to date, the well-posedness of multidimensional reflected BSDEs (or, systems of reflected BSDEs) has only been established in the case of convex reflection domains: see, e.g., {\color{black}\cite{Gegout-Petit-Pardoux-96, Klimsiak-Rozkosz-Slominski-15, Chassagneux-Richou-17,Fakhouri-Ouknine-Ren-18}}.
The systems of reflected BSDEs in convex domains appear in certain types of stochastic control problems, such as the switching problems: see, among others, {\color{black}\cite{Hamadene-Zhang-10,Hu-Tang-10,Chassagneux-Elie-Kharroubi-10,Chassagneux-Elie-Kharroubi-11,Martyr-16,Benezet-Chassagneux-Richou-19}}. On the other hand, in a certain class of control-stopping stochastic differential games, the associated equilibria are described by the systems of reflected BSDEs in non-convex domains, as shown, e.g., in \cite{GaydukNadtochiy3} (see also \cite{GaydukNadtochiy2} for the convex case). {\color{black} We also refer to \cite{Briand-Hibon-21}, which considers another case of a system of reflected BSDEs in a non-convex domain.}
This paper presents the first general well-posedness result for the systems of reflected BSDEs in non-convex domains under the assumption of a weak star-shape property, see Assumption \ref{ass:intro.domain} below.

In addition to the stochastic control-stopping games, the reflected BSDEs in non-convex domains have a direct connection to the theory of martingales on manifolds. We refer to \cite{Emery-89} for an introduction and an overview of this theory. One of the key questions therein is the following: given a random variable $\xi$ with values in a manifold $M$, is it possible to define a martingale $Y$ in $M$ such that the terminal value of this martingale (at time $T>0$) is given by $\xi$ (i.e., $Y_T=\xi$), and is such a martingale unique? A positive answer to this question, in particular, allows one to extend the notion of conditional expectation and gives one possible way to define a barycenter on a manifold (see e.g. \cite{Emery-Mokobodski-91,Picard-94}). We refer to \cite{Kendall-90,Kendall-91,Picard-91,Darling-95} for other applications. As explained in \cite{Darling-95}, it is possible to give a positive answer to this question by solving a BSDE with quadratic non-linearities with respect to the $z$-variable, stated in $\R^d$ -- the Euclidean space in which the manifold is embedded. It turns out that for a certain class of non-convex reflection domains $\cD$, the reflected BSDE in $\cD$ gives rise to a martingale on the manifold $\partial\cD$, see Section \ref{sect:sphere}. In particular, our results provide a new proof of the existence and uniqueness of a martingale with a prescribed terminal value in a given strict sub-sector of a hemisphere of $\mathbb{S}^{d-1}$, in the Markovian framework or under additional smallness assumptions (see the example in Section \ref{sect:sphere}).
 
On a technical level, our analysis is connected to the theory of BSDEs with quadratic growth in the $z$-variable. This connection is made precise in Section \ref{se:Markov}, but it can also be seen if one attempts to map a given non-convex domain into a convex one -- the resulting reflected BSDE in a convex domain will have quadratic terms in $z$. Thus, the reflected BSDEs in non-convex domains can be viewed as the quadratic reflected BSDEs in convex domains. This observation also explains the additional challenges of the case of a non-convex domain, relative to a convex one: the mathematical difficulties in the former case are similar to those arising in the well-posedness theory for the systems of quadratic BSDEs {\color{black} \cite{Tevzadze-08,Hu-Tang-16,Xing-Zitkovic-18,Harter-Richou-19}}. The present work uses some of the results developed in the latter theory: in particular, the results of \cite{Xing-Zitkovic-18} are crucial for our analysis.
 
 Another important connection is to the methods of \textcolor{black}{\cite{lions1981construction, Lions-Sznitman-84}}, which establish the well-posedness of the forward (or, regular) stochastic differential equations (SDEs) reflected at the boundary of a given domain. In particular, we use the arguments of \cite{Lions-Sznitman-84} to establish the stability of the solutions to the reflected BSDEs considered herein, see Section \ref{se:geometry}. It is important to mention, however, that many crucial arguments used in the proof of the well-posedness of a reflected (forward) SDE cannot be applied to the case of a reflected BSDE due to the adaptedness issues which, in particular, prohibit the application of the Skorokhod's mapping, used in \cite{Lions-Sznitman-84}, and of the standard localization methods.
 


The remainder of this paper is organized as follows. Section \ref{subse:setup} states the reflected BSDE (equation \eqref{eq reflected bsde}) and the main assumptions (Assumptions \ref{ass:intro.domain} and \ref{ass:intro.data}) which hold throughout the paper. Section \ref{se:geometry} describes various auxiliary properties and a priori estimates, as well as the stability (Proposition \ref{prop:stability2}) and uniqueness (Corollary \ref{th:uniqueness}) of the solutions to the target reflected BSDE in a certain class. Section \ref{se:Markov} describes a sequence of penalized quadratic BSDEs in a Markovian framework, shows that their solutions converge to a solution of the target reflected BSDE, and verifies that this solution belongs to the class in which the uniqueness holds, thus establishing the well-posedness of the target reflected BSDE in a Markovian framework (Theorem \ref{thm:existence:Markov}). In Section \ref{se:nonMarkov}, we approximate a general reflected BSDE by the Markovian ones, to obtain the well-posedness of the former (Theorem \ref{thm:general}) under an additional smallness assumption (Assumption \ref{ass:general}). Finally, Section \ref{sect:sphere} provides a more detailed description of the connection between the reflected BSDEs in non-convex domains and the martingales on manifolds, which, in particular, illustrates the sharpness of some of our assumptions.

\subsection{The setup and main assumptions}
\label{subse:setup}
Let $\cD$ be a subset of $\R^d$ given by
$$
\cD = \{y\in\R^d:\, \phi(y)<0\},
$$
with a function $\phi:\R^d\rightarrow\R$.
We denote by $\nabla$ the gradient, and by $\nabla^2$ the Hessian, of a given function.
For any subset $A$ of an Euclidean space, we denote its closure by $\bar A$ and, when $A \neq \emptyset$, we denote by $d(.,A)$ the distance function to $A$.


\begin{Assumption}\label{ass:intro.domain}
We assume that $\phi$ satisfies the following:
\begin{itemize}
\item (Compactness) There exists $R>0$, s.t. $\phi(y)>0$ for all $|y|\geq R$.
\item (Smoothness) $\phi\in C^2(\R^d)$, $|\nabla\phi(y)|>0$ for all $y\in\partial\cD$, and $\nabla^2\phi$ is locally Lipschitz.
\item (Weak star-shape property) There exists a non-empty open convex set $\mathcal{C} \subset \cD$ such that
\begin{itemize}
\item $0 \in \mathcal{C}$,
\item there exists a convex function $\phi_{\mathcal{C}}: \R^{d} \rightarrow \R$ satisfying $\phi_{\mathcal{C}} \in C^2(\R^{d})$,  
$$\mathcal{C} = \{y\in\R^d:\, \phi_{\mathcal{C}}(y)<0\},$$
$\phi_{\mathcal{C}}\geqslant \phi_{\mathcal{C}}(0)$ and $\phi_{\mathcal{C}}(y) = |y-\mathfrak{P}_{\bar{\mathcal{C}}}(y)|$ for all $y \in \R^d \setminus \mathcal{C}$ where $\mathfrak{P}_{\bar{\mathcal{C}}}$ stands for the projection function onto $\bar{\mathcal{C}}$,
\item it holds that
\begin{equation}
\label{def:gamma}
\gamma:= \inf_{y \in \partial\cD} \nabla \phi_{\mathcal{C}}(y)  \cdot \frac{\nabla \phi(y)}{|\nabla \phi(y)|}>0.
 \end{equation}

\end{itemize}
\end{itemize}
\end{Assumption}


\begin{Remark}\label{rem:weakStar.vs.strongStar}
If $\cD$ is a strictly star-shaped domain w.r.t. $0$, i.e., if it satisfies 
$$
\inf_{y \in \partial \cD} \frac{y}{|y|} \cdot \frac{\nabla \phi(y)}{|\nabla \phi(y)|}>0,
$$
then the weak star-shape property is also satisfied, with $\mathcal{C}$ being a ball of radius $\varepsilon>0$ centered at $0$, and with 
$$\phi_{\mathcal{C}}(y)=\varrho_\varepsilon(|y|-\varepsilon),$$
where $\varrho_\varepsilon:\mathbb{R} \rightarrow \mathbb{R}$ is a convex increasing function satisfying $\varrho_\varepsilon \in C^2(\mathbb{R},\mathbb{R})$, $\varrho_\varepsilon(x) = -\varepsilon/2$ for $x<-\varepsilon$ and  $\varrho_\varepsilon(x)=x$ for $x>0$.
\end{Remark}

\medskip

All stochastic processes and random variables, appearing in this paper, are constructed on a fixed stochastic basis $(\Omega,\mathbb{F},\P)$, with the filtration $\mathbb{F}$ being a completion of the natural filtration of a multidimensional Brownian motion $W$ in $\R^{d'}$ on a time interval $[0,T]$.

\smallskip

For $p\geq1$, we denote by $\mathcal{L}^p$ the space of (classes of equivalence of)\footnote{We drop this clarification in further definitions.} $\mathcal{F}_T$-measurable random variables $\xi$ (with values in a Euclidean space), s.t. $\|\xi\|_{\mathcal{L}^p}:=\esp{|\xi|^p}^{1/p}<\infty$. The space $\mathcal{L}^\infty$ stands for all $\mathcal{F}_T$-measurable essentially bounded random variables.
We also define $\HP{2}$ as the space of progressively measurable processes (with values in a Euclidean space) $Z$, s.t. $\|Z\|_{\HP{2}}:=\esp{ \int_0^T |Z_t|^2 \ud t }^{1/2}<\infty$.
Next, for $p\geq1$, we define $\mathcal{M}^p$ as the space of all continuous local martingales $M$ with $\|M\|_{\mathcal{M}^p}:=\esp{\langle M\rangle_T^{p/2}}^{1/p}<\infty$.
For $p \in [1,\infty]$, we denote by $\SP{p}$ the set of continuous adapted process $U$ such that $\NL{p}{\sup_{t \in[0,T]} |U_t|} < \infty$.
We also denote by $\mathrm{Var}_{t}(K)$ the variation of a process $K_\cdot$ (with values in a Euclidean space) on the time interval $[0,t]$ and by $\mathscr{K}^p$, for $p \in [1,\infty]$, the set of finite-variation process $K$ such that $\NL{p}{\mathrm{Var}_{[0,T]}(K)} < \infty$ and $K_0=0$.
Finally, we denote by $\BP{2}$ the set of processes $V \in \HP{2}$, satisfying
\begin{align*}
 \NB{2}{V} := \NL{\infty}{\mathrm{sup}_{t \in [0,T]} \esp{\int_{t}^T |V_s|^2 \ud s | \cF_{t}}}^\frac12<+\infty.  \;
\end{align*}
Let us remark that $V \in \BP{2}$ implies that the martingale $\int_0^. V_s \ud W_s$ is a BMO martingale, and $\NB{2}{V}$ is the BMO norm of $\int_0^. V_s \ud W_s$. We refer to \cite{Kazamaki-94} for further details about BMO martingales.

\medskip

We are investigating the well-posedness of the following reflected BSDE $(Y,Z,K) \in \SP{2}\times \HP{2} \times \mathscr{K}^1$ 
\begin{align}
\label{eq reflected bsde}
 \left\{ \begin{aligned}
 &(i)\;Y_t = \xi+\int_t^T f(s,Y_s,Z_s)\ud s-\int_t^T \ud K_s-\int_t^T Z_s \ud W_s, \quad 0 \leqslant t \leqslant T,\\
 & (ii)\; Y_t \in \bar{\cD} \text{ a.s.},\quad K_t=\int_0^t \mathfrak{n}(Y_s) \ud \mathrm{Var}_s(K),
  \quad 0 \leqslant t \leqslant T,
\end{aligned}
 \right.
\end{align}
where $\mathfrak{n}$ is the unit outward normal to $\partial\cD$, extended as zero into $\cD$:
$$
\mathfrak{n}(y) = \frac{\nabla\phi(y)}{|\nabla\phi(y)|},\quad \forall y\in\partial\cD\quad \textrm{and} \quad \mathfrak{n}(y)=0, \quad \forall y \in \cD.
$$

\smallskip

\begin{Assumption}\label{ass:intro.data}
We assume that $\xi$ takes values in $\bar\cD$, $f(\cdot,y,z)$ is progressively measurable, $f(t,\cdot,\cdot)$ is globally Lipschitz ($K_{f,y}$-Lipschitz in $y$ and $K_{f,z}$-Lipschitz in $z$), uniformly in $(t,\omega)$, and $\NL{\infty}{|f(\cdot,0,0)|}<\infty$. In addition, w.l.o.g. (in view of the boundedness of $\cD$), there exists a compact $\mathcal{K}\subset \R^d$, s.t. $f(t,y,z)=0$ whenever $y\notin\mathcal{K}$.
\end{Assumption}







\medskip

Assumptions \ref{ass:intro.domain} and \ref{ass:intro.data} hold throughout the rest of the paper even if not cited explicitly.

\section{Geometric properties and a priori estimates}
\label{se:geometry}

In this section, we derive useful geometric properties of the domain $\cD$, expressed via the corresponding properties of the function $\phi$, and construct an auxiliary function $\psi$ which is used in the next section to define a sequence of approximating equations to \eqref{eq reflected bsde}.
We also present some key \emph{a priori} estimates and properties of the solutions to the RBSDEs \eqref{eq reflected bsde}.

\subsection{Absolute continuity of the process $K$}
As noticed in \cite{Gegout-Petit-Pardoux-96}, we can take advantage of the smoothness of $\cD$ to show that the process $K$ is absolutely continuous with respect to the Lebesgue measure.

\begin{Lemma}\label{le abs cont K}
Assume that $(Y,Z,K) \in \SP{2} \times \HP{2} \times \mathscr{K}^1$ solves \eqref{eq reflected bsde}. Then, almost every path of
$K$ is absolutely continuous with respect to the Lebesgue measure.
\end{Lemma}
\proof
Applying It\^o's formula to $t \mapsto \phi(Y_t)$, we obtain
\begin{align}
&\ud \phi(Y_t) = \left(-\nabla \phi(Y_t)\cdot f(t,Y_t,Z_t) +\frac12\mathrm{Tr}[Z_t^\top\nabla^2\phi(Y_t)Z_t] \right) \ud t\nonumber\\
&+ \nabla \phi(Y_t)\cdot \ud K_t + \nabla \phi(Y_t)\cdot Z_t\ud W_t\label{eq ito phi}
\end{align}
Then, the It\^o-Tanaka formula applied to the positive part of the semi-martingale $-\phi(Y_t)$ reads
\begin{align}\label{eq ito-tanaka phi}
\ud [-\phi(Y_t)]^+ = \1_{\set{-\phi(Y_t)>0}} \ud [-\phi(Y_t)] + \frac12 \ud L^0_t,
\end{align}
where $L^0$ is the local time of the semi-martingale $-\phi(Y)$ at zero. Since $\phi(Y_t) \le 0$,  we have
$\ud [-\phi(Y_t)]^+ = -\ud \phi(Y_t)$ which yields, combining \eqref{eq ito phi}--\eqref{eq ito-tanaka phi},
\begin{align}
&\1_{\set{\phi(Y_t)=0}}\left(-\nabla \phi(Y_t)\cdot f(t,Y_t,Z_t) +\frac12\mathrm{Tr}[Z_t^\top\nabla^2\phi(Y_t)Z_t] \right) \ud t
+|\nabla \phi(Y_t)| \ud \mathrm{Var}_t(K) \nonumber
\\
& \quad +\1_{\set{\phi(Y_t)=0}}\nabla \phi(Y_t)\cdot Z_t\ud W_t + \frac12 \ud L^0_t =0 \;. \nonumber
\end{align}
In particular, we deduce that
\begin{align}
\label{est:dVarK}
|\nabla \phi(Y_t)| \ud \mathrm{Var}_t(K) \le \1_{\set{\phi(Y_t)=0}}\left[\nabla \phi(Y_t)\cdot f(t,Y_t,Z_t) -\frac12\mathrm{Tr}[Z_t^\top\nabla^2\phi(Y_t)Z_t] \right]^+ \ud t,
\end{align}
which proves the absolute continuity of $K$.
\eproof

\subsection{The exterior sphere property}

The following lemma states the well known observation that, for any boundary point of a smooth domain, there exists a small enough tangent external sphere, see e.g. \cite{Lions-Sznitman-84}.

\begin{Lemma}\label{le:intro.extSphere}
There exists $R_0>0$, s.t.
\begin{align}\label{eq lions-sznitman}
(y-y')\cdot \mathfrak{n}(y) + \frac1{2R_0}|y-y'|^2 \ge 0\;,\quad \forall\,y \in \partial \cD, \; y' \in \bar{\cD}.
\end{align}
\end{Lemma}
\begin{proof}
Due to the smoothness of $\phi$, for any $y\in\partial\cD$ and $y'\in\bar \cD$, there exists $\lambda\in[0,1]$, s.t.
\begin{equation}\label{eq.R0.pf.eq1}
0\geq \phi(y') = \phi(y) + (y'-y) \cdot \mathfrak{n}(y) |\nabla\phi(y)| + \frac{1}{2} (y - y')^\top \nabla^2\phi(\lambda y + (1-\lambda)y')
(y-y'),
\end{equation}
It only remains to notice that: $\phi=0$ and $|\nabla\phi|$ is bounded away from zero on $\partial\cD$, and $|\nabla^2\phi|$ is bounded from above on $\bar\cD$. Thus, we obtain the statement of the lemma.
\qed
\end{proof}

\smallskip

Using the above lemma, we can define the projection operator that is used in the subsequent sections. To this end, we first define the set
\begin{equation*}
Q=\{y\in\R^d:\, d(y,\cD) < R_0\},
\end{equation*}
and the set-valued projection operator
\begin{equation*}
\mathfrak{P}(y) = \text{argmin}_{x\in\bar\cD} |x-y|,\quad y\in\R^d.
\end{equation*}

\begin{Corollary}\label{cor:intro.K}
For any $y\in Q$, $\mathfrak{P}(y)$ is a singleton.
\end{Corollary}
\begin{proof}
It is easy to see that, for a ball $B_r(y)\subset\R^d$, with radius $r>0$ and center at $y$, we have:
\begin{equation}\label{eq.intro.cor1.eq1}
(x-x')\cdot \frac{y-x}{|y-x|} + \frac1{2r}|x-x'|^2 = 0,\quad \forall\,x,x'\in\partial B_r(y).
\end{equation}
Next, assume that there exist $y\in \R^d\setminus\bar\cD$ and $x\neq x'\in\bar\cD$, s.t.
$$
|x-y|=|x'-y| = \text{argmin}_{z\in\bar\cD} |z-y|.
$$
Then, it is clear that $x,x'\in \partial B_r\cap \partial\cD$, with $r=\min_{z\in\bar\cD} |z-y|<R_0$, and the equations \eqref{eq lions-sznitman}, \eqref{eq.intro.cor1.eq1} yield a contradiction.
\qed
\end{proof}

\smallskip

W.l.o.g., we will identify the value of $\mathfrak{P}(y)$ with its only element, for any $y\in Q$.

\smallskip

\begin{Remark}
 \label{rem:domination:VarK}
Thanks to \eqref{eq.R0.pf.eq1} and to \eqref{est:dVarK} we easily deduce for any solution $(Y,Z,K) \in \SP{2} \times \HP{2} \times \mathscr{K}^1$ to \eqref{eq reflected bsde} that
\begin{align*}
\ud \mathrm{Var}_t(K) \le \1_{\set{\phi(Y_t)=0}}\left(\left[\frac{\nabla \phi(Y_t)}{|\nabla \phi(Y_t)|}\cdot f(t,Y_t,Z_t)\right]^+ +\frac{1}{2R_0}|Z_t|^2 \right)\ud t,
\end{align*}
with $R_0$ satisfying \eqref{eq lions-sznitman}.
\end{Remark}

\smallskip

The following lemma provides an alternative to \eqref{eq reflected bsde}(ii), and it becomes useful in the subsequent sections.

\begin{Lemma}
\label{lem alternative characterisation}
Assume that $(Y,Z,K) \in \SP{2} \times \HP{2} \times \mathscr{K}^1$ solves \eqref{eq reflected bsde}(i) and $Y_t \in \bar{\cD}$ a.s. for all $t\in[0,T]$. Then 
$$
K_t = \int_0^t \mathfrak{n}(Y_s) \ud \mathrm{Var}_s(K),\quad t\in[0,T],
$$
holds if and only if there exists a constant $c>0$, depending only on $\cD$, such that for all essentially bounded continuous adapted process $V$ in $\cD$,  we have
\begin{align}
\int_0^T \left( (Y_s - V_s) +c|Y_s-V_s|^2 \mathfrak{n}(Y_t) \right) \ud K_s \ge0
\;.
\end{align}
\end{Lemma}
\proof
One implication is a direct consequence of Lemma \ref{le:intro.extSphere}. The other implication is a mere generalization of Lemma 2.1 in \cite{Gegout-Petit-Pardoux-96}.
\eproof



\subsection{The pseudo-distance function}
\label{subse:pseudo-distance}

In this subsection, we modify the function $\phi$ in order to construct a new smooth function $\psi$ which satisfies the inequality \eqref{def:gamma} in $\mathbb{R}^d \setminus \cD$ instead of $\partial \cD$.
We denote by $\vartheta:\R\rightarrow[0,1]$ an infinitely smooth nondecreasing function which is equal to zero on $(-\infty,0]$ and to one on $[1,\infty)$. We also choose a large enough $R>1$, s.t. $\cD\subset B_{R-1}(0)$, and a small enough $\epsilon\in(0,1)$, s.t., for all $y\in B_{R+1}(0)\setminus\cD$, we have:
$$
\phi(y)\leq\epsilon\quad\Rightarrow\quad y\in B_R(0),\quad \nabla \phi_{\mathcal{C}}(y)\cdot\nabla\phi(y) >0.$$ 
Then, we define
\begin{equation}\label{eq.intro.psi.def}
\tilde\phi(y) := \phi^+(y) (1-\vartheta(|y|-R-1)) + \vartheta(|y|-R),
\quad\psi(y) := \tilde\phi(y) + \kappa|y| \vartheta(\tilde\phi(y)/\epsilon),
\end{equation}
for an arbitrary constant $\kappa>0$.

\smallskip

We refer to $\psi$ as the pseudo-distance function.

\smallskip

Notice that
\begin{align*}
\nabla \phi_{\mathcal{C}}(y)\cdot\nabla\psi(y) &= \nabla \phi_{\mathcal{C}}(y)\cdot \nabla\tilde\phi(y) 
+ \kappa \nabla \phi_{\mathcal{C}}(y)\cdot \frac{y}{|y|} \vartheta(\tilde\phi(y)/\epsilon)\\
&+ \kappa \nabla \phi_{\mathcal{C}}(y)\cdot\nabla\tilde\phi(y) |y| \vartheta'(\tilde\phi(y)/\epsilon)/\epsilon\\
&= \nabla \phi_{\mathcal{C}}(y)\cdot \nabla\tilde\phi(y) \left(1 + \kappa|y| \vartheta'(\tilde\phi(y)/\epsilon)/\epsilon\right) + \kappa \nabla \phi_{\mathcal{C}}(y)\cdot\frac{y}{|y|} \vartheta(\tilde\phi(y)/\epsilon).
\end{align*}
It is clear that $\psi\in C^2(\R^d\setminus\bar\cD)$ and that its derivatives up to the second order are locally Lipschitz-continuous on $\R^d\setminus\cD$.
It is also easy to see that $\tilde\phi(y)\in(0,\epsilon]$ if and only if $y\in B_{R+1}(0)\setminus\cD$ and $\phi(y)\leq\epsilon$, in which case $y\in B_R(0)$, $\tilde{\phi}(y) = \phi(y)$, $\nabla\tilde{\phi}(y) = \nabla\phi(y)$, and
$$
\nabla \phi_{\mathcal{C}}(y)\cdot\nabla\psi(y) \geq \nabla \phi_{\mathcal{C}}(y)\cdot\nabla\phi(y) >0,
$$
where we also observed that $\inf_{y \in \R^d\setminus \cD} \nabla \phi_{\mathcal{C}}(y)\cdot y/|y|>0$, which follows from the convexity of $\mathcal{C}$ and from the fact that $0 \in \mathcal{C}$.
If $\tilde\phi(y)\leq0$, then $y\in\bar\cD$. If $\tilde\phi(y)>\epsilon$, then
$$
\nabla \phi_{\mathcal{C}}(y)\cdot\nabla\psi(y) = \nabla \phi_{\mathcal{C}}(y)\cdot \nabla\tilde\phi(y) +  \kappa \nabla \phi_{\mathcal{C}}(y)\cdot\frac{y}{|y|},
$$
which can be made positive for all $y\in\R^d\setminus \cD$ by choosing large enough $\kappa>0$, as $|\nabla\tilde\phi|$ is bounded on $\R^d\setminus \cD$ and $\inf_{y \in \R^d\setminus \cD} \nabla \phi_{\mathcal{C}}(y)\cdot y/|y|>0$.

The following lemma summarizes the above properties of $\psi$ and states several additional properties which can be easily verified.

\begin{Lemma}\label{le:intro.psi.prop}
There exist constants $R,\epsilon,\kappa>0$, s.t. the function $\psi$ defined in \eqref{eq.intro.psi.def} satisfies the following properties.
\begin{enumerate}
\item $\psi$ is globally Lipschitz-continuous in $\R^d$.
\item There exist constants $c,C>0$, s.t. $c\, d(y,\cD) \le \psi(y) \leq C\,d(y,\cD)$ for $y\in\R^d$.
\item $\psi\in C^2(\R^d\setminus\bar\cD)$, and its derivatives up to the second order are globally Lipschitz-continuous in $\R^d\setminus\cD$.
\item $\inf_{y\in\R^d\setminus\cD}\nabla \phi_{\mathcal{C}}(y)\cdot\nabla\psi(y)>0$.
\item $\inf_{y\in\R^d\setminus\cD} |\nabla \psi(y)| >0$.
\item $\psi(y)=\phi(y)=0$, $\nabla\psi(y)=\nabla\phi(y)$, and $\nabla^2\psi(y)= \nabla^2\phi(y)$, for $y\in\partial\cD$.
\end{enumerate}
\end{Lemma}

In the remainder of the paper, we fix $(R,\epsilon,\kappa)$ as in the above lemma and consider the associated pseudo-distance function $\psi$. For convenience, we also extend the vector-valued function $\mathfrak{n}$ to $\R^d$ as follows:
$$
\mathfrak{n}(y) = \frac{1}{|\nabla\psi(y)|}\nabla\psi(y)\,\bone_{\set{\R^d\setminus\cD}}(y).
$$



\subsubsection{Asymptotic convexity of the squared pseudo-distance}

Due to Lemma \ref{le:intro.psi.prop}, the Hessian of $\psi^2$, denoted $\nabla^2\psi^2$, is well defined in $\R^d\setminus\cD$ (it is extended to the boundary of the latter set by continuity). The following lemma shows that $\nabla^2\psi^2$, viewed as a bilinear form, becomes positive semidefinite close to $\cD$.

\begin{Lemma}\label{le:intro.conv}
There exists a constant $C>0$, s.t., for all $y\in \R^d\setminus \cD$ and $z\in\R^d$,
$$
z^\top \nabla^2\psi^2(y) z \ge  -C \psi(y) |z|^2.
$$
\end{Lemma}
\begin{proof}
Notice that, for $y\in\R^d\setminus\cD$ and $z\in\R^d$,
$$
\nabla^2\psi^2(y) = 2\nabla\psi(y) \nabla^\top\psi(y) + 2\psi(y) \nabla^2\psi(y),
$$
$$
z^\top \nabla^2\psi^2(y) z
= 2 (\nabla\psi(y)\cdot z)^2 + 2\psi(y) z^\top \nabla^2\psi(y) z\geq  2\psi(y) z^\top \nabla^2\psi(y) z.
$$
Using the fact that $\nabla^2\psi$ is bounded (cf. the third property in Lemma \ref{le:intro.psi.prop}) and the second property in Lemma \ref{le:intro.psi.prop}, we complete the proof.
\qed
\end{proof}


\subsection{\emph{A priori} estimates}

In this subsection, we prove \emph{a priori} estimates in the case of general terminal condition $\xi$ and generator $f$. 
We first introduce the appropriate ``smallness assumption''.

\begin{Assumption}
 \label{ass:general}
 We assume that at least one of the following four conditions is fulfilled with some $\theta \ge 1$:
\begin{enumerate}[(i)]
 \item $|\phi_{\mathcal{C}}^+(\xi)|_{\mathscr{L}^{\infty}}< \frac{\gamma R_0}{\theta}$  and
 $\nabla \phi_{\mathcal{C}}(y)\cdot f(s,y,z) \leqslant 0,\quad \forall s,y,z \in [0,T] \times \bar{\cD}\setminus \mathcal{C} \times \R^{d \times d'}$,
 \item or $\sup_{x \in\cD} \phi_{\mathcal{C}}^+(x) < \frac{\gamma R_0}{\theta}$,
 \item or $\cC$ is the Euclidean ball centered at $0$  with  radius $\lambda>0$, and 
 $$
 |\xi|^2_{\mathscr{L}^{\infty}} < \lambda^2 + \frac{2R^2_0}{\theta},
 \quad\nabla \phi_{\mathcal{C}}(y)\cdot f(s,y,z) \leqslant 0,
 \quad \forall s,y,z \in [0,T] \times \bar{\cD}\setminus \mathcal{C} \times \R^{d \times d'},
 $$
 \item or $\cC$ is the Euclidean ball centered at $0$  with  radius $\lambda>0$, and 
 $$
 \sup_{x \in \cD} |x|^2 < \lambda^2 + \frac{2R^2_0}{\theta}\,,
 $$
\end{enumerate}
with $R_0$ satisfying \eqref{eq lions-sznitman} and $\gamma$ appearing in Assumption \ref{ass:intro.domain}.
\end{Assumption}

It is worth mentioning that Assumption \ref{ass:general} is not our standing assumption and is cited explicitly whenever it is invoked. In particular, our well-posedness results in the Markovian framework do not require the smallness assumption, see Section \ref{se:Markov} .

\smallskip

Next, we consider the following class of solution:
\begin{Definition}\label{de uniqueness class}
For any $\theta\geq1$, we denote by $\mathscr{U}(\theta,\xi,f,T)$ the set of all solutions $(Y,Z,K) \in \SP{2}\times \HP{2} \times \mathscr{K}^1$ to \eqref{eq reflected bsde} s.t.
\begin{align}
  \label{de bound expo K 2}
   \mathbb{E}\left[ e^{\frac{\theta p}{R_0} \mathrm{Var}_T(K)} \right] < \infty,
\end{align}
with some $p>1$ and with $R_0$ satisfying \eqref{eq lions-sznitman}.
\end{Definition}

In the sequel, we will generally drop $(\xi,f,T)$ in the notation for the class $\mathscr{U}$. Note also that we will mainly consider $\theta = 1$ or $\theta = 2$.

\smallskip

The following proposition clarifies the link between Assumption \ref{ass:general} and the class $\mathscr{U}(\theta)$.

\begin{Proposition} \label{A priori estimate reflected BSDE}
 Let $(Y,Z,K) \in \SP{2}\times \HP{2}\times \mathscr{K}^1$ be a solution to the RBSDE \eqref{eq reflected bsde}. 
Then, $Z \in \BP{2}$.
 Moreover, if Assumption \ref{ass:general} holds for some $\theta \ge 1$, then, there exist constants $C$ and $p>1$, which depend only on $K_{f,y}$, $K_{f,z}$, $\gamma$, $\lambda$, $\sup_{y \in \cD} |y|$, $\sup_{y \in \bar{\cD}} \phi_{\mathcal{C}}^+(y)$, $\|\phi_{\mathcal{C}}^+(\xi)\|_{\mathscr{L}^{\infty}}$, $\|f(.,0,0)\|_{\mathscr{L}^{\infty}}$ and $R_0$ (recall Assumption \ref{ass:intro.domain} and \eqref{eq lions-sznitman}), such that 
\begin{equation}
  \label{bound expo K 2}
   \mathbb{E}\left[ e^{\frac{\theta p}{R_0} \mathrm{Var}_T(K)} \right] \leqslant C.
 \end{equation}
Thus, under Assumption \ref{ass:general}, any solution $(Y,Z,K)$ belongs to $\mathscr{U}(\theta,\xi,f,T)$.
\end{Proposition}
\proof
1. We start by applying It\^o-Tanaka's formula to $\phi_{\mathcal{C}}^+(Y_t)$ (note that $\phi_{\mathcal{C}}^+$ is convex): for all $t\leqslant t'$,
\begin{align}
 &\mathbb{E}_t\left[\int_t^{t'} \nabla \phi_{\mathcal{C}}^+(Y_s)\cdot \ud K_s\right] \leqslant  \mathbb{E}_t\left[ \phi_{\mathcal{C}}^+(Y_{t'})+\int_t^{t'} \nabla \phi_{\mathcal{C}}^+(Y_s)\cdot f(s,Y_s,Z_s) \ud s\right]. \label{ineq1}
\end{align}
In the equation above and the proofs below, we use the shorter notation $\EFp{t}{\cdot}$ for
$\esp{\cdot|\cF_t}$.
\\Recalling Assumption \ref{ass:intro.domain}, we also have
\begin{align} \label{ineq2} \gamma \mathbb{E}_t \left[ \int_t^{t'} \ud \mathrm{Var}_s (K)\right] \leqslant \mathbb{E}_t \left[ \int_t^{t'} \nabla \phi_{\mathcal{C}}^+(Y_s) \cdot \mathfrak{n}(Y_s)\ud \mathrm{Var}_s (K)\right]= \mathbb{E}_t \left[ \int_t^{t'} \nabla \phi_{\mathcal{C}}^+(Y_s)\cdot \ud K_s\right].
\end{align}
This yields, for $t'=T$, that
\begin{align}\label{eq main point for ass 5.2.i}
\gamma \mathbb{E}_t \left[ \int_t^{T} \ud \mathrm{Var}_s (K)\right] 
\le
 \mathbb{E}_t\left[ \phi_{\mathcal{C}}^+(\xi)+\int_t^{T} \nabla \phi_{\mathcal{C}}^+(Y_s)\cdot f(s,Y_s,Z_s) \ud s\right] \,.
\end{align}
Next, we consider an arbitrary $\varepsilon >0$ and apply It\^o's formula to $\varepsilon |Y_t|^2$ between $t$ and $t'$, to obtain:
\begin{align*}
 \varepsilon \mathbb{E}_t\left[\int_t^{t'} |Z_s|^2\ud s \right] \leqslant& \varepsilon \mathbb{E}_t\left[|Y_{t'}|^2 \right]+C_1\varepsilon \mathbb{E}_t\left[ \int_t^{t'}  (1+|f(s,0,0)|+|Z_s|)\ud s\right] \\
 & + \varepsilon C_1\mathbb{E}\left[\int_t^{t'} \ud \mathrm{Var}_s (K)\right],
\end{align*}
where we used that $|Y|$ is bounded. The above inequality implies
\begin{align}
 \frac{\varepsilon}{2}\mathbb{E}_t\left[\int_t^{t'} |Z_s|^2\ud s \right] \leqslant& \varepsilon \mathbb{E}_t\left[|Y_{t'}|^2 \right]+C_{\varepsilon} \mathbb{E}_t\left[ \int_t^{t'}  (1+|f(s,0,0)|)\ud s\right]+ \varepsilon C_1\mathbb{E}\left[\int_t^{t'} \ud \mathrm{Var}_s (K)\right]. \label{ineq3}
\end{align}
Setting $t'=T$ and $\varepsilon=1$ in the previous inequality and combining it with \eqref{eq main point for ass 5.2.i}, we obtain
\begin{align}
 \frac{1}{2}\mathbb{E}_t\left[\int_t^{T} |Z_s|^2\ud s \right] \leqslant&
   \mathbb{E}_t\left[|\xi|^2+   \frac{C_1}{\gamma}\phi_{\mathcal{C}}^+(\xi) \right]+C \mathbb{E}_t\left[ \int_t^{T}  (1+|f(s,0,0)|)\ud s\right] \nonumber
  \\
  &+  \frac{C}{\gamma}\mathbb{E}\left[\int_t^{T} \nabla \phi_{\mathcal{C}}^+(Y_s)\cdot f(s,Y_s,Z_s) \ud s\right]. \label{eq step for z bmo}
\end{align}
Now, we observe that
\begin{align*}
\frac{C_1}{\gamma}\mathbb{E}\left[\int_t^{T} \nabla \phi_{\mathcal{C}}^+(Y_s)\cdot f(s,Y_s,Z_s) \ud s\right]
\le C\mathbb{E}_t\left[ \int_t^{t'}  (1+|f(s,0,0)|)\ud s\right] + \frac14\mathbb{E}_t\left[\int_t^{T} |Z_s|^2\ud s \right] 
\end{align*}
Inserting the previous estimate back into \eqref{eq step for z bmo}, we get
\begin{align*}
 \frac{1}{4}\mathbb{E}_t\left[\int_t^{T} |Z_s|^2\ud s \right] \leqslant&
   \mathbb{E}_t\left[|\xi|^2+   \frac{C}{\gamma}\phi_{\mathcal{C}}^+(\xi) \right]+C \mathbb{E}_t\left[ \int_t^{T}  (1+|f(s,0,0)|)\ud s\right]\;.
\end{align*}
This proves that $Z \in \BP{2}$.

\smallskip

\noindent 2. We now turn to the estimation of the exponential moments of $\mathrm{Var}_T(K)$, under the smallness Assumption \ref{ass:general}.

\smallskip

\noindent 2.a First, combining Assumption \ref{ass:general}(i) with \eqref{eq main point for ass 5.2.i}, we obtain
\begin{align}
 \frac{\theta}{R_0}\left\|\sup_{t \in [0,T]}\mathbb{E}_t \left[ \int_t^{T} \ud \mathrm{Var}_s (K)\right]\right\|_{\mathscr{L}^\infty}  <  1.
\end{align}
Then, we apply the energy inequalities for non-decreasing processes with bounded potential (see, e.g., (105.1)-(105.2) in \cite{Dellacherie-Meyer-85}) to obtain \eqref{bound expo K 2} in this case. 

\smallskip

\noindent 2.b Let now Assumption \ref{ass:general}-(ii) hold. Using \eqref{ineq1}-\eqref{ineq2} and recalling that $|Y|$ is bounded, we obtain, for all $0 \leqslant t < t' \leqslant T$ and for any $\varepsilon>0$,
\begin{align*}
\gamma  \mathbb{E}_t \left[ \int_t^{t'} \ud \mathrm{Var}_s (K)\right] \leqslant  \sup_{y \in \bar{\cD}} \phi_{\mathcal{C}}^+(y)+C_{\varepsilon}'(t'-t)(1+|f(.,0,0)|_{\mathscr{L}^{\infty}})+\frac{\varepsilon}{2} \mathbb{E}_t \left[ \int_t^{t'} |Z_s|^2 \ud s\right]. 
\end{align*}
Using the above inequality and \eqref{ineq3} (with the same $\varepsilon>0$), we obtain: 
\begin{align*}
 (\gamma-C\varepsilon)  \mathbb{E}_t \left[ \int_t^{t'} \ud \mathrm{Var}_s (K)\right] \leqslant \varepsilon \sup_{y \in \bar{\cD}} |y|^2 + \sup_{y \in \bar{\cD}} \phi_{\mathcal{C}}^+(y)+C_{\varepsilon}''(t'-t)(1+|f(.,0,0)|_{\mathscr{L}^{\infty}}).
\end{align*}
In particular, by taking $\varepsilon$ small enough, we conclude that, for any $\varepsilon'>0$, there exists $C_{\varepsilon'}>0$ such that
\begin{align}
\mathbb{E}_t \left[ \int_t^{t'} \ud \mathrm{Var}_s (K)\right] \leqslant  \frac{\sup_{y \in \bar{\cD}} \phi_{\mathcal{C}}^+(y)}{\gamma}(1+\varepsilon')+C_{\varepsilon'}(t'-t).\label{eq.nonMarkov.mainThm.case2Est}
\end{align}
\smallskip

Next, using  \eqref{eq.nonMarkov.mainThm.case2Est} and Assumption \ref{ass:general}(ii), we conclude that there exist $0 <\varepsilon''<1$, $p>1$, and $N \geqslant 1$, depending only on $K_{f,y}$, $K_{f,z}$, $\gamma$, $\sup_{y \in \cD} |y|$, $\sup_{y \in \bar{\cD}} \phi_{\mathcal{C}}(y)^+$,  $\|f(.,0,0)\|_{\mathscr{L}^{\infty}}$ and $R_0$, such that, a.s.:
\begin{align}
\mathbb{E}_t \left[ \int_t^{T(k+1)/N} \ud \mathrm{Var}_s (K)\right] \leqslant \frac{R_0}{\theta p}(1-\varepsilon''), \quad \forall 0\leqslant k <N, \,\, \forall t \in [Tk/N,T(k+1)/N].\label{eq.nonMarkov.mainThm.K.small}
\end{align}
Then, we apply the energy inequalities for non-decreasing processes with bounded potential (see, e.g., (105.1)-(105.2) in \cite{Dellacherie-Meyer-85}), to obtain, for all $0 \leqslant k <N$,
\begin{align}\label{eq pre-iterate}
\mathbb{E}_{Tk/N} \left[e^{\frac{\theta p}{R_0}\int_{Tk/N}^{T(k+1)/N} \ud \mathrm{Var}_s (K)}\right] \leqslant \tilde{C},
\end{align}
with $\tilde{C}$ that depends only on $K_{f,y}$, $K_{f,z}$, $\gamma$, $\sup_{y \in \cD} |y|$, $\sup_{y \in \bar{\cD}} \phi_{\mathcal{C}}(y)^+$,  $\|f(.,0,0)\|_{\mathscr{L}^{\infty}}$ and $R_0$. 
We now observe that
\begin{align*}
\esp{e^{\frac{\theta p}{R_0}{Var}_T(K)}} &= \esp{e^{\frac{ \theta p}{R_0} \mathrm{Var}_{T(N-1)/N} (K)} 
\mathbb{E}_{T(N-1)/N} \left[e^{\frac{ \theta p}{R_0}\int_{T(N-1)/N}^{T} \ud \mathrm{Var}_s (K)}\right] 
} 
\\
&\le \tilde{C} \esp{e^{\frac{\theta p}{R_0} \mathrm{Var}_{T(N-1)/N} (K)} },
\end{align*}
where we used \eqref{eq pre-iterate} with $k=N-1$ to obtain the last inequality.
Iterating the above procedure concludes the proof of this case.

\smallskip

\noindent 2.c 
Next, we let Assumption \ref{ass:general}(iii) hold. Using \eqref{est:dVarK}, the linear growth of $f$, and Young's inequality, we have, for all $\varepsilon>0$,
\begin{align}
\label{caseiii:est:varT}
 \mathrm{Var}_T(K) \le C_{\varepsilon}+\frac{1+\varepsilon}{2R_0}\int_0^T\1_{\set{\phi(Y_t)=0}}|Z_t|^2 \ud t.
\end{align}
Moreover, we apply It\^o-Tanaka formula to $(|Y_s|^2-\lambda^2)^+$ to obtain, for all $t\leqslant t'$,
\begin{align}
  \nonumber
 &\mathbb{E}_t\left[\int_t^{t'} \1_{\set{\phi(Y_s)=0}}|Z_s|^2 \ud s \right]  \leqslant  \mathbb{E}_t\left[ \int_t^{t'} \1_{\set{\phi_{\cC}(Y_s)>0}}|Z_s|^2 \ud s \right]\\
 \leqslant &\mathbb{E}_t\left[ (|Y_{t'}|^2-\lambda^2)^+ + 2\int_t^{t'} \1_{\set{\phi_{\cC}(Y_s)>0}} |Y_s|\nabla\phi_{\cC} (Y_s)\cdot f(s,Y_s,Z_s) \ud s  \right], \label{caseiii:estimZ}
\end{align}
where we also recall that $\1_{\set{\phi_{\cC}(Y_s)>0}}|Y_s|\nabla\phi_{\cC} (Y_s)=\1_{\set{\phi_{\cC}(Y_s)>0}}Y_s$ since $\cC$ is a Euclidean ball centered at zero.
Then, by taking $t'=T$ in \eqref{caseiii:estimZ} and using Assumption \ref{ass:general}(iii), we obtain, for $\varepsilon>0$ and $p>1$ small enough, 
\begin{align*}
\frac{\theta p(1+\varepsilon)}{2R^2_0}\left\|\sup_{t \in [0,T]}\mathbb{E}_t \left[ \int_t^{T} \1_{\set{\phi(Y_s)=0}}|Z_s|^2 \ud s \right]\right\|_{\mathscr{L}^\infty}  <  1.
\end{align*}
It remains to apply the John-Nirenberg inequality for BMO Martingales (see Theorem 2.2 in \cite{Kazamaki-94}) and recall \eqref{caseiii:est:varT}, to conclude that 
$$\mathbb{E}\left[ e^{\frac{\theta p}{R_0} \mathrm{Var}_T(K)} \right] \leqslant
C_{\varepsilon}\mathbb{E}\left[ e^{\frac{\theta p (1+\varepsilon)}{2R^2_0}  \int_0^T \1_{\set{\phi(Y_s)=0}}|Z_s|^2 \ud s}\right] <+\infty,$$
which yields \eqref{bound expo K 2}.

\smallskip

\noindent 2.d Finally, the proof in the case of Assumption \ref{ass:general}(iv) follows from \eqref{caseiii:est:varT} and \eqref{caseiii:estimZ}, by partitioning $[0,T]$ into small time intervals as in step 2.b. To avoid the redundant calculations, we skip the details.
%
\eproof

%

\subsection{Stability and uniqueness in $\mathscr{U}(\theta)$}

Using the a priori estimates established in the previous subsection, we prove the following stability property of the solutions to \eqref{eq reflected bsde}.

\begin{Proposition}
\label{prop:stability2}
 Let us consider $(Y,Z,K) \in \SP{2}\times \HP{2}\times \mathscr{K}^1$ (resp. $(Y',Z',K')\in \SP{2}\times \HP{2}\times \mathscr{K}^1$) which solve the RBSDE \eqref{eq reflected bsde} with a domain $\cD$ (resp. $\cD'$), with a terminal condition $\xi$ (resp. $\xi'$), and with a generator $f$ (resp. $f'$). Assume, moreover, that there exists $p>1$ such that  
 \begin{equation}\label{eq main ass stab}
 \kappa :=  \mathbb{E}\left[e^{\frac{p}{R_0}(\mathrm{Var}_T (K)+\mathrm{Var}_T (K'))}\right] <+\infty,
 \end{equation}
with $R_0$ satisfying \eqref{eq lions-sznitman} for $\cD$ and $\cD'$. 
Let us denote by $\bar{\mathfrak{P}}$ (resp. $\bar{\mathfrak{P}}'$) a measurable selection of the projection operator onto $\cD$ (resp. $\cD'$). Then, the following stability result holds: there exists a constant $C>0$, which depends only on $K_{f,y}$, $K_{f',y}$, $K_{f,z}$, $K_{f',z}$ (recall Assumption \ref{ass:intro.data}), $\sup_{y \in \cD \cup \cD'} |y|$, $R_0$, and on $\kappa$,
and is such that
\begin{align*}
& \|Y-Y'\|_{\SP{2}}+\|Z-Z'\|_{\HP{2}}+\|K-K'\|_{\SP{2}}\\ \leqslant &C\mathbb{E}[|\xi-\xi'|^{2p/(p-1)}]^{(p-1)/(2p)} + C\mathbb{E}\left[ \left(\int_0^T |f(s,Y_s,Z_s)-f'(s,Y_s,Z_s)|\ud s\right)^{2p/(p-1)} \right]^{(p-1)/(2p)}\\
&+C\mathbb{E}\left[\sup_{s \in [0,T]}|Y_s-\bar{\mathfrak{P}}'(Y_s)|^{p/(p-1)}\right]^{(p-1)/(2p)}+C\mathbb{E}\left[\sup_{s \in [0,T]}|Y'_s-\bar{\mathfrak{P}}(Y'_s)|^{p/(p-1)}\right]^{(p-1)/(2p)}.
\end{align*}
\end{Proposition}
\proof
We apply It\^o's formula to the process
$$e^{\beta t+\frac{1}{R_0}(\mathrm{Var}_t(K)+\mathrm{Var}_t(K'))}|Y_t - Y'_t|^2,$$
with the constant $\beta$ to be determined later on. 
By denoting
$$
\df_t := f(t,Y_t,Z_t) - f'(t,Y'_t,Z'_t),\quad \delta \xi := \xi -\xi',
$$ 
$$
\Gamma_t := e^{\beta t+\frac{1}{R_0}(\mathrm{Var}_t(K)+\mathrm{Var}_t(K'))},\quad \dY := Y - Y',\quad \dZ = Z - Z',
$$
we obtain
\begin{align}
 &\Gamma_t|\dY_t|^2 + \int_t^T\!\! \Gamma_s |\dZ_s|^2 \ud s\nonumber \\
 =& \Gamma_T |\delta \xi |^2+ 2 \int_t^T\!\! \Gamma_s \dY_s\cdot \df_s \ud s - 2\int_t^T \Gamma_s \dY_s \cdot \ud K_s + 2\int_t^T \!\! \Gamma_s \dY_s \cdot \ud K'_s \nonumber \\
& -\beta \int_t^T \!\! \Gamma_s |\dY_s|^2 \ud s -\frac{1}{R_0} \int_t^T \!\! \Gamma_s |\dY_s|^2 \ud\mathrm{Var}_s (K) -\frac{1}{R_0} \int_t^T \!\! \Gamma_s |\dY_s|^2 \ud\mathrm{Var}_s(K')  \nonumber \\
&- 2 \int_t^T \!\! \Gamma_s\dY_s\cdot \dZ_s \ud W_s. \label{ito:stabilite}
\end{align}

Using the BDG inequality, the fact that $|\dY|$ is bounded, H\"older inequality (with $q=p/(p-1)>1$ being the conjugate exponent), we obtain:
\begin{align*}
\esp{\sup_{t \in [0,T]} \left|\int_0^t \!\! \Gamma_s\dY_s\cdot \dZ_s \ud W_s\right|}
&\le C \esp{\left(\int_0^T  \left|\Gamma_s \dZ_s \right|^2 \ud s\right )^\frac12}
\\
&\le C \esp{(\Gamma_T)^p}^{\frac1p} \esp{ (\int_0^T |\dZ_s|^2 \ud s)^{\lceil q/2 \rceil}}^\frac1q<\infty,
\end{align*}
where the last inequality is due to \eqref{eq main ass stab} and to the Energy Inequality (since $Z,Z' \in \BP{2}$). Hence, we conclude that the local martingale term in the right hand side of \eqref{ito:stabilite} is a true martingale. 

Next, we estimate the second term in the right hand side of \eqref{ito:stabilite} using the Lipschitz property of $f'$:
\begin{align*}
 \dY_s\cdot \df_s  &\leqslant |\delta Y_s||f(s,Y_s,Z_s)-f'(s,Y_s,Z_s)|+\beta |\dY_s|^2+\frac{1}{4}|\dZ_s|^2,
\end{align*}
provided $\beta>0$ is large enough.
In addition, the condition \eqref{eq reflected bsde}(ii) and the exterior sphere property (recall \eqref{eq lions-sznitman}) yield
\begin{align*}
&- 2\int_t^T \Gamma_s \dY_s \cdot \ud K_s  -\frac{1}{R_0} \int_t^T \!\! \Gamma_s |\dY_s|^2 \ud\mathrm{Var}_s (K)\\
= & - 2\int_t^T \Gamma_s (\bar{\mathfrak{P}}(Y'_s) - Y'_s) \cdot \ud K_s 
- 2\int_t^T \Gamma_s\left(Y_s - \bar{\mathfrak{P}}(Y'_s)\right) \cdot \ud K_s
-\frac{1}{R_0} \int_t^T \!\! \Gamma_s |\dY_s|^2 \ud\mathrm{Var}_s (K)\\
= & - 2\int_t^T \Gamma_s\left(Y_s - \bar{\mathfrak{P}}(Y'_s)\right) \cdot \ud K_s
-\frac{1}{R_0} \int_t^T \!\! \Gamma_s |Y_s - \bar{\mathfrak{P}}(Y'_s)|^2 \ud\mathrm{Var}_s (K)\\
& + \frac{1}{R_0} \int_t^T \!\! \Gamma_s (|\bar{\mathfrak{P}}(Y'_s) - Y'_s|^2 - |\dY_s|^2) \ud\mathrm{Var}_s (K)
- 2\int_t^T \Gamma_s (\bar{\mathfrak{P}}(Y'_s) - Y'_s) \cdot \ud K_s\\
\leqslant &  C\int_t^T \!\! \Gamma_s | \bar{\mathfrak{P}}(Y'_s)-Y'_s| \ud\mathrm{Var}_s (K)
\leqslant  C \Gamma_T \sup_{s \in [0,T]} | \bar{\mathfrak{P}}(Y'_s)-Y'_s| ,
\end{align*}
where for the last inequality, we used:
$
\int_t^T\exp\left(\frac{\mathrm{Var}_s (K)}{R_0}\right) \ud\mathrm{Var}_s (K) \le R_0 \exp\left(\frac{\mathrm{Var}_T (K)}{R_0}\right)\,.
$
By the same arguments we obtain
\begin{align*}
 2\int_t^T \!\! \Gamma_s \dY_s \cdot \ud K'_s -\frac{1}{R_0} \int_t^T \!\! \Gamma_s |\dY_s|^2 \ud\mathrm{Var}_s(K'))& \leqslant C \Gamma_T \sup_{s \in [0,T]} | \bar{\mathfrak{P}}'(Y_s)-Y_s| .
\end{align*}
Using the above estimates, we take expectations on both sides of \ref{ito:stabilite}, with $t=0$, and apply H\"older inequality to obtain 
\begin{align}
\|\Gamma^{1/2}\delta Z\|_{\HP{2}} &\leqslant \mathbb{E}\left[\Gamma_T |\delta \xi|^2 +2 \int_0^T \Gamma_s|\delta Y_s||f(s,Y_s,Z_s)-f'(s,Y_s,Z_s)|\ud s\right]^{1/2}\nonumber \\
&\quad +\mathbb{E}\left[\Gamma_T \left(\sup_{s \in [0,T]} | \bar{\mathfrak{P}}'(Y_s)-Y_s|+| \bar{\mathfrak{P}}(Y'_s)-Y'_s|\right) \right]^{1/2} \nonumber \\
&\leqslant C\mathbb{E}[|\delta \xi|^{2q}]^{1/(2q)} \nonumber \\
&+ 2\mathbb{E}\left[\sup_{s \in [0,T]}(\Gamma_s^{1/2}|\delta Y_s|)\Gamma_T^{1/2}\int_0^T |f(s,Y_s,Z_s)-f'(s,Y_s,Z_s)|\ud s \right]^{1/2} \nonumber\\
&\quad +C\mathbb{E}\left[ \sup_{s \in [0,T]} \left(|\bar{\mathfrak{P}}'(Y_s)-Y_s|+| \bar{\mathfrak{P}}(Y'_s)-Y'_s|\right)^q\right]^{1/(2q)}. \label{stab Z 2}
\end{align}
Using \eqref{ito:stabilite} and \eqref{stab Z 2}, we apply BDG, H\"older and Young inequalities to obtain
\begin{align}
\|\Gamma^{1/2} \delta Y\|_{\SP{2}} 
&\leqslant C\mathbb{E}[|\delta \xi|^q]^{1/(2q)}\nonumber \\
&+ C\mathbb{E}\left[ \sup_{s \in [0,T]}(\Gamma_s^{1/2}|\delta Y_s|)\Gamma_T^{1/2}\int_0^T |f(s,Y_s,Z_s)-f'(s,Y_s,Z_s)|\ud s \right]^{1/2} \nonumber\\
&\quad +C\mathbb{E}\left[ \sup_{s \in [0,T]} \left(|\bar{\mathfrak{P}}'(Y_s)-Y_s|+| \bar{\mathfrak{P}}(Y'_s)-Y'_s|\right)^q\right]^{1/(2q)}\nonumber\\
\leqslant C\mathbb{E}[|\delta \xi|^q]^{1/(2q)}& + \frac{1}{2}\|\Gamma^{1/2} \delta Y\|_{\SP{2}}
 + C\mathbb{E}\left[ \left(\int_0^T |f(s,Y_s,Z_s)-f'(s,Y_s,Z_s)|\ud s\right)^{2q} \right]^{1/(2q)}\nonumber \\
&\quad + C\mathbb{E}\left[ \sup_{s \in [0,T]} \left(|\bar{\mathfrak{P}}'(Y_s)-Y_s|+| \bar{\mathfrak{P}}(Y'_s)-Y'_s|\right)^q\right]^{1/(2q)}.
\label{stab Y 2}
\end{align}
Then, combining \eqref{stab Z 2}, Young inequality, and \eqref{stab Y 2}, yields
\begin{align}
 &\|Y-Y'\|_{\SP{2}}+\|Z-Z'\|_{\HP{2}}
 \leqslant \|\Gamma^{1/2}\delta Y\|_{\SP{2}}+\|\Gamma^{1/2}\delta Z\|_{\HP{2}} \nonumber \\
 \leqslant &C\mathbb{E}[|\xi-\xi'|^{2q}]^{1/(2q)} + C\mathbb{E}\left[ \left(\int_0^T |f(s,Y_s,Z_s)-f'(s,Y_s,Z_s)|\ud s\right)^{2q} \right]^{1/(2q)}\nonumber \\
 & +C\mathbb{E}\left[ \sup_{s \in [0,T]} \left(|\bar{\mathfrak{P}}'(Y_s)-Y_s|+| \bar{\mathfrak{P}}(Y'_s)-Y'_s|\right)^q\right]^{1/(2q)}.\label{stab Y Z}
\end{align}

Finally, we recall that
$$K_t - K'_t = \dY_t-\dY_0 + \int_0^t f(s,Y_s,Z_s) - f'(s,Y_s',Z'_s)\ud s-\int_0^t \dZ_s \ud W_s.$$
Then, the BDG inequality, the Lipschitz property of $f'$, as well as \eqref{stab Y Z}, yield
\begin{align*}
\|K-K'\|_{\SP{2}} \leqslant&  C\mathbb{E}[|\xi-\xi'|^{2q}]^{1/(2q)} + C\mathbb{E}\left[ \left(\int_0^T |f(s,Y_s,Z_s)-f'(s,Y_s,Z_s)|\ud s\right)^{2q} \right]^{1/(2q)}\\
& +C\mathbb{E}\left[ \sup_{s \in [0,T]} \left(|\bar{\mathfrak{P}}'(Y_s)-Y_s|+| \bar{\mathfrak{P}}(Y'_s)-Y'_s|\right)^q\right]^{1/(2q)},
\end{align*}
which completes the proof of the proposition.
\eproof

\smallskip

In a general non-Markovian framework, we obtain the following uniqueness result as a direct consequence of Proposition \ref{prop:stability2}.
\begin{Corollary}
\label{th:uniqueness}
 The reflected BSDE \eqref{eq reflected bsde} has at most one solution  in the class $\mathscr{U}(2)$.
\end{Corollary}
\proof Indeed, it suffices to check that, for any two solutions in the class $\mathscr{U}(2)$, \eqref{eq main ass stab} holds.
This follows directly from the Cauchy-Schwarz inequality.
\eproof

This uniqueness result is improved in the Markovian setting: see Theorem \ref{thm:existence:Markov} and Remark \ref{rem:unique.Markov}.



\section{Well-posedness in a Markovian framework}
\label{se:Markov}

In this section, we establish the existence and uniqueness of the solution to \eqref{eq de penalised equation} under the assumption that the terminal condition and the generator of the reflected BSDE are functions of a Markov diffusion process $X$ in $\R^{d'}$:
\begin{equation}
\label{SDE}
 X_t=x + \int_0^t b(s, X_s) \ud s + \int_0^t \sigma(s,X_s) \ud W_s,\quad x\in\R^{d'}.
\end{equation}

Namely, we make the following assumptions.

\begin{Assumption}
\label{ass:intro:SDE}
We assume that $(b,\sigma)$ are bounded measurable functions, uniformly Lipschitz with respect to $x$, and such that $\sigma^\top\sigma$ is uniformly positive definite (i.e., uniformly elliptic), which implies in particular that $\sigma$ is invertible.
\end{Assumption}

\begin{Assumption}
 \label{ass:markov}
We assume that
$$
 \xi :=g(X_T) \quad \text{and} \quad f(t,y,z):=F(t,X_t,y,z),
$$
where $g$ is $\alpha$-H\"older and $\bar{\cD}$-valued, $F$ is measurable in all variables, globally Lipschitz in $(y,z)$, and s.t. $|F(\cdot,\cdot,0,0)|$ is bounded.
\end{Assumption}

Recall that Assumptions \ref{ass:intro.domain} and \ref{ass:intro.data} hold throughout the paper, even if they are not cited explicitly.

\subsection{Penalized equation}
We begin by noticing that $\psi^2\in C^1(\R^d)$ and denote 
$$
\Psi(y):=\frac{1}{2}\nabla\psi(y)^2 = \psi(y) \nabla\psi(y),\quad y\in\R^d,
$$
where we extend (naturally) $\nabla\psi$ to $\cD$ by zero. We also extend $\nabla^2\psi^2$ to $\cD$ by zero.

\smallskip

It is useful to note that there exist constants $c,C$, s.t.
\begin{equation}\label{eq.pen.Psi.psi}
0<c \psi \leq |\Psi| \leq C\psi.
\end{equation}

\smallskip

Next, we consider the following penalized equation:
\begin{align}
\label{eq de penalised equation}
Y^n_t = \xi + \int_{t}^Tf(s,Y^n_s,Z^n_s) \ud s-  \int_t^T n \Psi(Y^n_s)(1 + |Z^n_s|^2)  \ud s - \int_t^T Z^n_s \ud W_s.
\end{align}
Let us remark that, contrarily to the convex framework tackled in \cite{Gegout-Petit-Pardoux-96}, it is natural (and necessary) to add a $|z|^2$ inside the penalization term due to \eqref{est:dVarK}.
For convenience, we introduce the notation:
\begin{align}
\Phi^{n}_t := \int_0^t n \Psi(Y^n_s) \ud s,\quad&\Theta^{n}_t := \int_0^t n \Psi(Y^n_s) |Z^n_s|^2 \ud s,\nonumber\\
&K^n_t := \Phi^{n}_t + \Theta^{n}_t.\label{eq.Kn.def}
\end{align}

\subsection{Existence of a solution to the penalized equation}

We start by considering the following family of approximating BSDEs, indexed by a pair of positive integers $M=(M_1,M_2)$:
\begin{align}
&Y_t^{n,M}= g(X_T) + \int_t^T F^{n,M}(s,X_s,Y_s^{n,M},Z_s^{n,M})\ud s-\int_t^T Z_s^{n,M}\ud W_s,\label{eq.penalized.trunc}
\end{align}
with 
$$
F^{n,M}(t,x,y,z) := f(t,x,y,z) - n\rho_{M_1}(\psi(y))\nabla\psi(y)(1+\rho_{M_2}(|z|^2)),
\quad\rho_k(x):=x \wedge k.
$$ 
The above BSDE has a globally Lipschitz generator and, therefore, is known to have a unique Markovian solution $(Y^{n,M},Z^{n,M})\in \SP{2}\times \HP{2}$ (see, e.g., Theorem 4.1 in \cite{ElKaroui-Peng-Quenez-97}).
The following Proposition uses the weak star-shape property of $\cD$, stated in Assumption \ref{ass:intro.domain}, to establish a uniform estimate on $(Y^{n,M},Z^{n,M})$.
 
 \begin{Lemma}\label{pr first estimates}
 There exists a constant $C>0$, 
 s.t., for any $n\geq1$, any $M=(M_1,M_2)$, and any $t\in[0,T]$, the following holds a.s.:
 \begin{align}
  |Y^{n,M}_t|^2 + \EFp{t}{\int_t^T |Y^{n,M}_s|^2 + |Z^{n,M}_s|^2 \ud s} \le 
  C\, \EFp{t}{|\xi|^2+\int_t^T (1+|f(s,0,0)|^2) \ud s},\label{eq.aprioriEst1.eq2}
 \end{align}
 \begin{align}
 \EFp{t}{\int_t^T n \rho_{M_1}(\psi(Y^{n,M}_s)) (1+\rho_{M_2}(|Z_s^{n,M}|^2)) \ud s}
 \le C\, \EFp{t}{|\xi|^2+\int_t^T (1+ |f(s,0,0)|^2) \ud s}.\label{eq.aprioriEst1.eq3}
 \end{align} 
 \end{Lemma}
 
 
 \proof
W.l.o.g., we assume that $\phi_{\mathcal{C}}$ attains its minimum at zero.
Then, we consider arbitrary $t \in [0, T]$ and constants $\alpha>0$, $\beta > 0$, to be fixed later, and define
$$
 [t,T]\times \R^d \ni (s,y) \mapsto h(s,y):=e^{\beta(s-t)}\left(\alpha|y|^2+(\phi_{\mathcal{C}}(y) - \phi_{\mathcal{C}}(0))^2\right) \in \R.
$$
We observe that $(\phi_{\mathcal{C}} - \phi_{\mathcal{C}}(0))^2$ is convex and $h(s,y) \le e^{\beta(T-s)}c_0 |y|^2$, for some positive constant $c_0$.
Then, we apply It\^o's formula to the process $h(s,Y_s^{n,M})$ (recalling \eqref{eq de penalised equation}), to obtain
\begin{align}
 \alpha |Y^{n,M}_t|^2 &\leqslant h(t,Y^{n,M}_t)\leqslant h(T,\xi)\nonumber \\
 & + 2\int_t^T e^{\beta(s-t)}(\alpha Y^{n,M}_s+(\phi_{\mathcal{C}}(Y_s^{n,M})-\phi_{\mathcal{C}}(0))\nabla \phi_{\mathcal{C}}(Y_s^{n,M}))\cdot f(s,Y^{n,M}_s,Z^{n,M}_s) \ud s\nonumber\\
&-   \int_t^T n  \rho_{M_1}(\psi(Y^{n,M}_s)) \nabla_y h(s,Y_s^{n,M})\cdot\ \nabla\psi(Y^{n,M}_s) (1+ \rho_{M_2}(|Z^{n,M}_s|^2))  \ud s\label{eq.penal.est1.eq1}\\
&- 2\int_t^T \nabla_y h(s,Y_s^{n,M})\cdot Z^{n,M}_s \ud W_s 
- \alpha \int_t^T e^{\beta(s-t)}|Z^{n,M}_s|^2  \ud s -\beta \int_t^T e^{\beta(s-t)} |Y^{n,M}_s|^2 \ud s.\nonumber
\end{align}

\smallskip

As $Y^{n,M}\in\SP{2}$ and $Z^{n,M}\in\HP{2}$, the local martingale in the above representation is in $\mathcal{M}^1$ and, hence, is a true martingale.

\smallskip

Next, we notice that the fourth property in Lemma \ref{le:intro.psi.prop} implies the existence of a constant $c_1>0$, s.t.
$$
\nabla \phi_{\mathcal{C}}(Y^{n,M}_s)\cdot \nabla\psi(Y^{n,M}_s) \geq c_1 \1_{\set{Y_s^{n,m} \notin \cD}}.
$$
Then, there exist constants $c_2,c_3>0$ such that
\begin{align}
& -   \int_t^T n  \rho_{M_1}(\psi(Y^{n,M}_s)) \nabla_y h(s,Y_s^{n,M})\cdot\ \nabla\psi(Y^{n,M}_s) (1+ \rho_{M_2}(|Z^{n,M}_s|^2))  \ud s\nonumber\\
\leqslant&  -   2\int_t^T n e^{\beta(s-t)} \rho_{M_1}(\psi(Y^{n,M}_s)) (c_1 (\phi_{\mathcal{C}}(Y_s^{n,M})-\phi_{\mathcal{C}}(0))-\alpha c_2 |Y_s^{n,M}|) (1+ \rho_{M_2}(|Z^{n,M}_s|^2))  \ud s\nonumber\\
\leqslant&  -   2\int_t^T n e^{\beta(s-t)} \rho_{M_1}(\psi(Y^{n,M}_s)) \left[c_1 (\phi_{\mathcal{C}}(Y_s^{n,M})-\phi_{\mathcal{C}}(0))\right.
\label{eq.penalized.firstEst.weakStar.eq1}\\
&\left.\phantom{????????????}-\alpha c_2 (\phi_{\mathcal{C}}(Y_s^{n,M})+|\mathfrak{P}_{\bar{\mathcal{C}}}(Y_s^{n,M})|)\right] \left(1+ \rho_{M_2}(|Z^{n,M}_s|^2)\right)  \ud s\nonumber\\
\leqslant& -c_3 \int_t^T n e^{\beta(s-t)} \rho_{M_1}(\psi(Y^{n,M}_s))  (1+ \rho_{M_2}(|Z^{n,M}_s|^2))  \ud s,\nonumber
\end{align}
provided $\alpha$ is small enough. In the rest of the proof, we assume that $\alpha$ is chosen so that the above inequality holds.

\smallskip

Next, we remark that
\begin{align*}
&\left| 2\int_t^T e^{\beta(s-t)}\left(\alpha Y^{n,M}_s+(\phi_{\mathcal{C}}(Y_s^{n,M})-\phi_{\mathcal{C}}(0)\right)\nabla \phi_{\mathcal{C}}(Y_s^{n,M}))\cdot f(s,Y^{n,M}_s,Z^{n,M}_s) \ud s \right|\\
&\leq C_1\int_t^T e^{\beta(s-t)} \left((\alpha+1)|Y^{n,M}_s|-\phi_{\mathcal{C}}(0)\right)\left( |f(s,0,0)| + C_2 |Y^{n,M}_s| + C_2 |Z^{n,M}_s| \right) \ud s\\
&\leq \int_t^T e^{\beta(s-t)}\left(C_3|Y^{n,M}_s|^2 +C_3+ |f(s,0,0)|^2 + \frac{\alpha}{2}|Z^{n,M}_s|^2 \right) \ud s.
\end{align*}

\smallskip

\smallskip
Combining the above estimates and \eqref{eq.penal.est1.eq1}, we conclude that, for a large enough $\beta>0$, there exists a constant $C_4>0$, s.t. 
\begin{align*}
&\alpha |Y^{n,M}_t|^2 + \EFp{t}{c_3 \int_t^T n \rho_{M_1}(\psi(Y^{n,M}_s)) \left(1 + \rho_{M_2}(|Z^{n,M}_s|^2)\right)\ud s + \frac{\alpha}{2}\int_t^T |Y^{n,M}_s|^2 + |Z^{n,M}_s|^2 \ud s}\\
&\phantom{???????????????????????????}\le e^{\beta (T-t)}\EFp{t}{c_0|\xi|^2+\int_t^T (C_4+|f(s,0,0)|^2) \ud s},
\end{align*}
which yields the statement of the lemma.
\eproof

\medskip

\begin{Proposition}\label{prop:penal.WellPosed}
Under Assumptions \ref{ass:intro:SDE} and \ref{ass:markov}, for any $n\geq1$, the BSDE \eqref{eq de penalised equation} has a Markovian solution $(Y^n,Z^n)$. In particular, there exists a measurable function $u^n$ such that $Y^n_t=u^n(t,X_t)$. Moreover, the estimates \eqref{eq.aprioriEst1.eq2}--\eqref{eq.aprioriEst1.eq3} hold with $(Y^{n,M},Z^{n,M})$ and $\rho_{M_i}$ replaced, respectively, by any solution $(Y^{n},Z^{n})$ of \eqref{eq de penalised equation} and by the identity function.
\end{Proposition}
\proof
The main statement of the proposition follows from Theorem 2.8 in \cite{Xing-Zitkovic-18} (without the localization). To be able to apply the latter theorem, we first consider the following auxiliary BSDE, which can be viewed as a middle ground between \eqref{eq de penalised equation} and \eqref{eq.penalized.trunc}:
\begin{align}
&\tilde Y_t^{n,M_1}= g(X_T) + \int_t^T \tilde F^{n,M_1}(s,X_s,\tilde Y_s^{n,M_1},\tilde Z_s^{n,M_1})\ud s-\int_t^T \tilde Z_s^{n,M_1}\ud W_s,\label{eq.penal.midBSDE}
\end{align}
with 
$$
\tilde F^{n,M_1}(t,x,y,z) := f(t,x,y,z)-n\rho_{M_1}(\psi(y))\nabla\psi(y)(1+|z|^2)
$$
and $\rho_{M_1}(x)=x \wedge M_1$.
We claim that the unique solution $(Y^{n,M},Z^{n,M})$ of \eqref{eq.penalized.trunc} converges (along a subsequence) to a Markovian solution $(\tilde Y^{n,M_1},\tilde Z^{n,M_1})$ of \eqref{eq.penal.midBSDE}, as $M_2\rightarrow\infty$. Indeed, this claim follows directly from Theorem 2.8 in \cite{Xing-Zitkovic-18}. To verify the assumptions of the latter theorem, we first notice that, due to \eqref{eq.aprioriEst1.eq2}, there exists a constant $c>0$ such that $|Y^{n,M}_t|\leqslant c$, for all $t \in [0,T]$ and $n,M$. Moreover, for large enough $C>0$ (independent of $n$ and $M$), $h(y):=C\left(\alpha |y|^2+(\phi_{\mathcal{C}}(y)-\phi_{\mathcal{C}}(0))^2\right)$ is a global $c$-Lyapunov function for $(F^{n,M})_{M}$, in the sense of Definition 2.3 in \cite{Xing-Zitkovic-18}, where $\alpha$ is the constant chosen in the proof of Lemma \ref{pr first estimates}. Indeed, there exists a large enough $C>0$, s.t., for all $|y| \leqslant c$, we have:
\begin{align*}
&\frac{1}{2}C\, \mathrm{Tr}[(z\sigma)^\top(\nabla^2 h(y)) z\sigma] - C \nabla h(y) \cdot F^{n,M}\\
&\geqslant C\alpha\, \mathrm{Tr}[(z\sigma)^\top z\sigma] - 2 C \left[\alpha y+(\phi_{\mathcal{C}}(y)-\phi_{\mathcal{C}}(0))\nabla \phi_{\mathcal{C}}(y)\right]\cdot f(t,x,y,z)\\
&\quad+ 2Cn \left[\alpha y+(\phi_{\mathcal{C}}(y)-\phi_{\mathcal{C}}(0))\nabla \phi_{\mathcal{C}}(y)\right] \cdot\nabla\psi(y) \rho_{M_1}(\psi(y))(1+\rho_{M_2}(|z|^2))\\
&\geq |z|^2-C',
\end{align*}
where we used the uniform ellipticity of $\sigma^\top\sigma$, Assumption \ref{ass:intro.data}, and the fourth property in Lemma \ref{le:intro.psi.prop}, and repeated the estimates used in \eqref{eq.penalized.firstEst.weakStar.eq1}.
In addition, we have $|F^{n,M}(t,x,y,z)| \leqslant C+C_n |z|^2$, with the constants $(C,C_n)$ independent of $M_2$. Observing that $F^{n,M}$ converges to $\tilde F^{n,M_1}$ locally uniformly, as $M_2\rightarrow\infty$, we conclude that the assumptions of Theorem 2.8 in \cite{Xing-Zitkovic-18} are satisfied and that \eqref{eq.penal.midBSDE} has a Markovian solution $(\tilde Y^{n,M_1},\tilde Z^{n,M_1})$ which is a limit point of $\{(Y^{n,M},Z^{n,M})\}_{M_2}$.

Next, we recall that, due to \eqref{eq.aprioriEst1.eq2}, $|Y^{n,M}|$ is bounded uniformly over $M$. Hence, $|\tilde Y^{n,M_1}|$ can be bounded uniformly over $M_1\geq 1$, and, in turn, $(\tilde Y^{n,M_1},\tilde Z^{n,M_1})$ solve \eqref{eq de penalised equation} for any large enough $M_1>0$.

The estimates \eqref{eq.aprioriEst1.eq2}--\eqref{eq.aprioriEst1.eq3} are obtained by repeating the proof of Lemma \ref{pr first estimates} for the equation \eqref{eq de penalised equation} in place of \eqref{eq.penalized.trunc}.
\eproof

\subsection{A priori estimates}
\label{subse:aprioriestimates}

The following result relies on the asymptotic convexity of the squared pseudo-distance function, stated in Lemma \ref{le:intro.conv}.

\begin{Lemma} \label{pr penal pseudo distance}
Under Assumptions \ref{ass:intro:SDE} and \ref{ass:markov}, there exists a constant $C>0$, s.t., for any $n\geq1$, any solution $(Y^n,Z^n)$ of \eqref{eq de penalised equation}, and any $t\in[0,T]$, the following holds a.s.:
\begin{align*}
n\psi^2(Y^n_t)
+ &\EFp{t}{\int_t^T n^2 |\Psi(Y^n_s)|^2 \left( 1 + |Z^n_s|^2\right) \ud s}\\
&\phantom{?????????????????????????}\le C\, \EFp{t}{|\xi|^2 + \int_t^T|f(s,0,0)|^2 \ud s},
\end{align*}
and, in particular,
\begin{align*}
d(Y^n_t,\cD)\le C n^{-1/2}.
\end{align*}
\end{Lemma}
\proof
We begin by applying It\^o's formula to $|\psi(Y^n_t)|^2$, to obtain
\begin{align}
\psi^2(Y^n_t) &= 2 \int_{t}^T \Psi(Y^n_s)\cdot f(s,Y^n_s,Z^n_s) \ud s - 2 \int_t^T n |\Psi(Y^n_s)|^2 \ud s
\nonumber\\
&-  2\int_t^T n |\Psi(Y^n_s)|^2  |Z^n_s|^2  \ud s - 2 \int_t^T \Psi(Y^n_s)\cdot Z^n_s \ud W_s
\label{eq.pen.est2.0.eq1}\\
&- \frac{1}{2} \int_t^T \mathrm{Tr}[(Z^n_s)^\top \nabla^2\psi^2(Y^n_s) Z^n_s] \ud s\nonumber
\end{align}

\smallskip

\begin{Remark}\label{rem:explain.Ito}
Note that the Hessian of $\psi^2$ has a discontinuity at $\partial\cD$. To justify the use of It\^o's formula, we approximate $\psi^2$ by a sequence of $C^2$ functions $\{g^m\}$, so that $g^m$, $\nabla g^m$ and $\nabla^2 g^m$ converge, respectively, to $\psi^2$, $\nabla\psi^2$ and $\nabla^2\psi^2$ everywhere in $\R^d$, and $|\nabla g^m|$, $|\nabla^2 g^m|$ are locally bounded uniformly over $m$.
To construct such a sequence, we first define
\begin{equation*}
\hat\phi(y) := \phi(y) (1-\vartheta(|y|-R-1)) + \vartheta(|y|-R),
\quad\hat\psi(y) := \hat\phi(y) + \kappa |y| \vartheta(\hat\phi(y)/\epsilon),
\quad y\in\R^d,
\end{equation*}
where we recall the original function $\phi$, appearing in Assumption \ref{ass:intro.domain}, and use the same $\vartheta$, $R$, and $\epsilon$, as the ones used in Subsection \ref{subse:pseudo-distance} to define $\psi$ (see \eqref{eq.intro.psi.def}). It is clear that $\hat\psi(y)=\psi(y)$, for $y\in\R^d\setminus\cD$, and that $\hat\psi(y)=\phi(y)$, for $y\in\cD$. Thus, $\hat\psi$ is a smooth extension of $\psi$ into $\cD$.
Next, we consider an infinitely smooth nondecreasing function $\rho:\R\rightarrow\R$, such that $\rho(x)=-1$ for $x\leq-1$ and $\rho(x)=x$ for $x\geq0$, and define
$$
g^m(y):=\frac{1}{m^2} \rho^2\left(m\hat\psi(y)\right),\quad y\in\R^d.
$$
It is easy to check by a direct computation that $g^m(y)$, $\nabla g^m(y)$ and $\nabla^2 g^m(y)$ converge to zero as $m\rightarrow\infty$, for any $y\in\cD$. On the other hand, $g^m(y)$ and its first two derivatives coincide with $\psi^2(y)$ and with its respective derivatives, for all $y\in\R^d\setminus\cD$ and all $m$. Thus, we obtain the desired sequence $\{g^m\}$. Applying It\^o's formula to $g^m(Y^n_t)$ and using the dominated convergence theorem to pass to the limit as $m\rightarrow\infty$, we establish \eqref{eq.pen.est2.0.eq1}.
\end{Remark}

\smallskip

As $|\Psi|$ is linearly bounded (see Lemma \ref{le:intro.psi.prop}), we conclude, as in the proof of Lemma \ref{pr first estimates}, that the local martingale in the above representation is a true martingale.

\smallskip

Next, we note that
\begin{equation*}
2\Psi(Y^n_s)\cdot f(s,Y^n_s,Z^n_s)
\leq n|\Psi(Y^n_s)|^2 + n^{-1} |f(s,Y^n_s,Z^n_s)|^2,
\end{equation*}
and use Lemma \ref{pr first estimates}, to obtain:
\begin{align}
&\EE_t \int_{t}^T 2\Psi(Y^n_s)\cdot f(s,Y^n_s,Z^n_s) \ud s\label{eq.pen.est2.0.eq6}\\
&\phantom{????????????????????}
\leq \EE_t \int_{t}^T n|\Psi(Y^n_s)|^2 \ud s + C n^{-1} \EFp{t}{|\xi|^2+\int_t^T(1+|f(s,0,0)|^2) \ud s}.\nonumber
\end{align}

\smallskip

In addition, Lemmas \ref{le:intro.psi.prop} and \ref{le:intro.conv} yield
\begin{equation*}
\mathrm{Tr}[(Z^n_s)^\top \nabla^2\psi^2(Y^n_s) Z^n_s]
\geq - C \Psi(Y^n_s) |Z^n_s|^2.
\end{equation*}
Then,
\begin{align}
&-n\left(\frac{1}{2}\mathrm{Tr}[(Z^n_s)^\top \nabla^2\psi^2(Y^n_s) Z^n_s] + 2 n |\Psi(Y^n_s)|^2  |Z^n_s|^2\right)
\nonumber\\
&\phantom{?????????????????}\le 
\left(Cn|\Psi(Y^n_s)| -  c n^{2}|\Psi(Y^n_s)|^2 \right)|Z^n_s|^2\label{eq.pen.est2.0.eq2}
\end{align}
Next, we observe that
\begin{align}
&\left(Cn|\Psi(Y^n_s)| -  C n^{2}|\Psi(Y^n_s)|^2 \right)|Z^n_s|^2 \leq C |Z^n_s|^2.
\label{eq.pen.est2.0.eq3}
\end{align}
Collecting \eqref{eq.pen.est2.0.eq2}--\eqref{eq.pen.est2.0.eq3} and using Lemma \ref{pr first estimates}, we obtain
\begin{align}
&-\EE_t \int_t^T\left(\mathrm{Tr}[(Z^n_s)^\top \nabla^2\psi^2(Y^n_s) Z^n_s] + 2 n|\Psi(Y^n_s)|^2  |Z^n_s|^2\right) \ud s
\label{eq.pen.est2.0.eq5}\\
& \leq C n^{-1} \EE_t \int_t^T |Z^n_s|^2 \ud s
\leq C n^{-1} \EFp{t}{|\xi|^2+\int_t^T(1+|f(s,0,0)|^2) \ud s}.\nonumber
\end{align}
Taking the conditional expectation in \eqref{eq.pen.est2.0.eq1}, multiplying both sides by $n$, and using \eqref{eq.pen.est2.0.eq6}, \eqref{eq.pen.est2.0.eq5}, we complete the proof.
\eproof

\medskip

The following proposition improves the rate of convergence of $Y^n$ to $\cD$.

\begin{Proposition} \label{pr penal pseudo distance 2}
Under Assumptions \ref{ass:intro:SDE} and \ref{ass:markov}, there exist $\mathfrak{N},C>0$, s.t. for any $n\geq\mathfrak{N}$, any solution $(Y^n,Z^n)$ of \eqref{eq de penalised equation}, and any $t\in[0,T]$, the following holds a.s.:
\begin{align*}
n\psi(Y^n_t)\le C
\end{align*}
\end{Proposition}


\proof
First, we denote by $\| \nabla^2\psi(y)\|_*$ the maximum absolute value across all negative parts of the entries of the matrix $\nabla^2\psi(y)$.
Next, we fix arbitrary $\epsilon,\varepsilon>0$ satisfying
$$
\epsilon \leq \left(\varepsilon + \sup_{y\in \partial\cD}\frac{\|\nabla^2\psi(y)\|_*}{|\nabla \psi(y)|^2} \right)^{-1},
\quad \epsilon \leq  \left(\varepsilon+2\sup_{y\in\partial\cD,\,z\in\R^{d\times m},\,s\in[0,T]} \frac{\NL{\infty}{\mathfrak{n}(y)\cdot f(s,y,z)}}{|\nabla \psi(y)|}\right)^{-1},
$$
and define
$$
\Psi^n(y) := (\psi(y)-1/(\epsilon n))^+ \nabla\psi(y) = \frac{1}{2}\nabla \left((\psi(y)-1/(\epsilon n))^+\right)^2,
$$
$$
\tilde H(y):= \nabla^2 \left((\psi(y)-1/(\epsilon n))^+\right)^2
= 2\nabla\psi(y) \nabla^\top\psi(y)\,\bone_{\psi\geq 1/(\epsilon n)} + 2 (\psi(y)-1/(\epsilon n))^+ \nabla^2\psi(y).
$$

\smallskip

Next, we apply It\^o's formula to $((\psi(Y^n_t)-1/(\epsilon n))^+)^2$ (the validity of It\^o's formula for the function $((\psi-1/(\epsilon n))^+)^2$ is justified similarly to Remark \ref{rem:explain.Ito}), to obtain
\begin{align}
((\psi(Y^n_t)-1/(\epsilon n))^+)^2 &= 2 \int_{t}^T \Psi^n(Y^n_s)\cdot f(s,Y^n_s,Z^n_s) \ud s - 2 \int_t^T n |\Psi(Y^n_s)| |\Psi^n(Y^n_s)| \ud s
\nonumber\\
&-  2\int_t^T n |\Psi(Y^n_s)| |\Psi^n(Y^n_s)|  |Z^n_s|^2  \ud s - 2 \int_t^T \Psi^n(Y^n_s)\cdot Z^n_s \ud W_s
\label{eq.pen.est2.eq1}\\
&- \frac{1}{2} \int_t^T \mathrm{Tr}[(Z^n_s)^\top \tilde H(Y^n_s) Z^n_s] \ud s.\nonumber
\end{align}

\smallskip

As $|\Psi^n|$ is linearly bounded (see Lemma \ref{le:intro.psi.prop}), we conclude, as in the proof of Proposition \ref{pr first estimates}, that the local martingale in the above representation is a true martingale.

\smallskip

Next, we note that
\begin{align*}
&2\Psi^n(Y^n_s)\cdot f(s,Y^n_s,Z^n_s) - n |\Psi(Y^n_s)| |\Psi^n(Y^n_s)|\\
&\leq |\Psi^n(Y^n_s)| \left( 2|\mathfrak{n}(Y^n_s)\cdot f(s,Y^n_s,Z^n_s)| - n \psi(Y^n_s) |\nabla \psi(Y^n_s)|\right).\nonumber
\end{align*}
Notice that, whenever $\Psi^n(Y^n_s)>0$, we have $\psi(Y^n_s)\geq 1/(\epsilon n)$ and, hence, $n \psi(Y^n_s) \geq 1/\epsilon$. Then, since $\epsilon>0$ satisfies 
$$
\epsilon \leq  \left(\varepsilon+2\sup_{y\in\partial\cD,\,z\in\R^{d\times m},\,s\in[0,T]} \frac{\NL{\infty}{\mathfrak{n}(y)\cdot f(s,y,z)}}{|\nabla \psi(y)|}\right)^{-1},
$$
and since $Y^n_s$ is close to $\cD$ for large enough $n$ (due to Lemma \ref{pr penal pseudo distance}),
we conclude that
\begin{align*}
n \psi(Y^n_s) |\nabla \psi(Y^n_s)|
&\geq |\nabla \psi(Y^n_s)|
\left(\varepsilon + 2\sup_{y\in\partial\cD,\,z\in\R^{d\times m},\,s\in[0,T]} \frac{\NL{\infty}{\mathfrak{n}(y)\cdot f(s,y,z)}}{|\nabla \psi(y)|}\right)\\
&\geq |\nabla \psi(Y^n_s)| \varepsilon/2 + 2 |\mathfrak{n}(Y^n_s)\cdot f(s,Y^n_s,Z^n_s)|
\end{align*}
and, in turn,
\begin{align}\label{eq.pen.est2.eq2}
&2\Psi^n(Y^n_s)\cdot f(s,Y^n_s,Z^n_s) - n |\Psi(Y^n_s)| |\Psi^n(Y^n_s)|\leq0.
\end{align}

\smallskip

Next, we recall that
$$
\frac12\mathrm{Tr}[(Z^n_s)^\top \tilde H(Y^n_s) Z^n_s]
\geq \left(\psi(Y^n_s)-1/(\epsilon n)\right)^+ \mathrm{Tr}[(Z^n_s)^\top \nabla^2\psi(Y^n_s) Z^n_s]
$$
and, hence,
\begin{equation*}
-\frac12\mathrm{Tr}[(Z^n_s)^\top \tilde H(Y^n_s) Z^n_s]
\leq \|\nabla^2\psi(Y^n_s)\|_* \left(\psi(Y^n_s)-1/(\epsilon n)\right)^+ |Z^n_s|^2.
\end{equation*}
In addition,
$$
-|\Psi(Y^n_s)| |\Psi^n(Y^n_s)| |Z^n_s|^2
= -|\nabla \psi(Y^n_s)|^2 \psi(Y^n_s) \left(\psi(Y^n_s)-1/(\epsilon n)\right)^+ |Z^n_s|^2.
$$
Collecting the two equations above, we deduce
\begin{align*}
&-\frac12\mathrm{Tr}[(Z^n_s)^\top \tilde H(Y^n_s) Z^n_s] - n |\Psi(Y^n_s)| |\Psi^n(Y^n_s)|  |Z^n_s|^2
\nonumber\\
&\phantom{?????????????????}\le 
\left(\psi(Y^n_s)-1/(\epsilon n)\right)^+ |Z^n_s|^2
\left(\|\nabla^2\psi(Y^n_s)\|_* - |\nabla \psi(Y^n_s)|^2 n \psi(Y^n_s) \right).
\end{align*}
Recall that, whenever $\psi(Y^n_s)\geq 1/(\epsilon n)$, we have $n \psi(Y^n_s) \geq 1/\epsilon$.
Then, since $\epsilon>0$ satisfies 
$$
\epsilon \leq \left(\varepsilon+\sup_{y\in \partial\cD}\frac{\|\nabla^2\psi(y)\|_*}{|\nabla \psi(y)|^2} \right)^{-1},
$$
and since $Y^n_s$ is close to $\cD$ for large enough $n$,
we conclude that
\begin{align*}
|\nabla \psi(Y^n_s)|^2 n \psi(Y^n_s)
&\geq |\nabla \psi(Y^n_s)|^2 \left(\varepsilon+\sup_{y\in \partial\cD}\frac{\|\nabla^2\psi(y)\|_*}{|\nabla \psi(y)|^2} \right)\\
&\geq |\nabla \psi(Y^n_s)|^2 \varepsilon/2 + \|\nabla^2\psi(Y^n_s)\|_*
\end{align*}
and, in turn,
\begin{align}
&-\frac{1}{2} \mathrm{Tr}[(Z^n_s)^\top \tilde H(Y^n_s) Z^n_s] 
- n |\Psi(Y^n_s)| |\Psi^n(Y^n_s)|  |Z^n_s|^2\leq0.\label{eq.pen.est2.eq3}
\end{align}


\smallskip

Taking the conditional expectation in \eqref{eq.pen.est2.eq1}, we make use of equations \eqref{eq.pen.est2.eq2} and \eqref{eq.pen.est2.eq3}, and of the fact that $|\Psi^n|\leq |\Psi|$, to obtain
\medskip
\begin{align*}
((\psi(Y^n_t)-1/(\epsilon n))^+)^2
+ &\EE_t\int_t^T n |\Psi^n(Y^n_s)|^2 \left(1 + |Z^n_s|^2\right) \ud s\leq 0
\end{align*}
and complete the proof.
\eproof

\medskip

Using Proposition \ref{pr penal pseudo distance 2}, we can improve the statement of Proposition \ref{prop:penal.WellPosed} and deduce that the H\"older norm of the Markovian solution of the penalized BSDE is bounded uniformly over $n$.

\begin{Corollary}\label{cor:unifHolder}
Under Assumptions \ref{ass:intro:SDE} and \ref{ass:markov}, there exist constants $\mathfrak{N}\geq1$, $\alpha' \in (0,1]$, and $C>0$ (independent of $n$), s.t., for any $n\geq\mathfrak{N}$, the BSDE \eqref{eq de penalised equation} has a Markovian solution $(Y^n,Z^n)$, with $Y^n_t=u^n(t,X_t)$, and any such solution satisfies
\begin{align}\label{eq unif holder}
\sup_{(t,x)\neq(t',x')} \frac{|u^n(t,x)-u^n(t',x')|}{|t-t'|^{\alpha'/2}+|x-x'|^{\alpha'}} \leqslant C.
\end{align}
\end{Corollary}
\proof
The statement of the corollary follows from Theorem 2.5 in \cite{Xing-Zitkovic-18} (without the localization). To verify the assumptions of the latter theorem, we consider the following capped version of \eqref{eq de penalised equation}:
\begin{align}
&\hat Y_t^{n,N}= g(X_T) + \int_t^T \hat F^{n,N}(s,X_s,\hat Y_s^{n,N},\hat Z_s^{n,N})\ud s-\int_t^T \hat Z_s^{n,N}\ud W_s,\label{eq.penal.midBSDE.2}
\end{align}
with 
$$
\hat F^{n,N}(t,x,y,z) := f(t,x,y,z)-\rho_N(n\psi(y))\nabla\psi(y)(1+|z|^2)
$$
and $\rho_N(x)=x \wedge N$. 
Propositions \ref{prop:penal.WellPosed} and \ref{pr penal pseudo distance 2} imply the existence of (large enough) $N,\mathfrak{N}>0$, s.t., for every $n\geq\mathfrak{N}$, there exists a Markovian solution $(Y^{n},Z^{n})$ of \eqref{eq de penalised equation}, with $Y^n_t=u^n(t,X_t)$, and any such solution also solves \eqref{eq.penal.midBSDE.2}. Moreover, there exists $c>0$ such that $|u^n| \leqslant c$ for all $n$.

Next, we fix $N$ as above and verify easily (as in the proof of Proposition \ref{prop:penal.WellPosed}) that, for large enough $C>0$ and small enough $\alpha>0$ (independent of $n$), $C(\alpha |y|^2+ (\phi_{\mathcal{C}}(y)-\phi_{\mathcal{C}}(0))^2)$ is a global $c$-Lyapunov function for $(\hat F^{n,N})_{n}$, in the sense of Definition 2.3 in \cite{Xing-Zitkovic-18}. 
In addition, $|\hat F^{n,N}(t,x,y,z)| \leqslant C+C_N |z|^2$, with the constants $(C,C_N)$ independent of $n$.
Thus, Theorem 2.5 in \cite{Xing-Zitkovic-18} yields the uniform boundedness of the H\"older norm of $u^n$.
\eproof

\medskip

W.l.o.g. we assume that the statements of Proposition \ref{pr penal pseudo distance 2} and Corollary \ref{cor:unifHolder} hold with $\mathfrak{N}=1$. 
From Corollary \ref{cor:unifHolder}, we deduce that there exists a subsequence of $\{u^n\}_{n \ge 1}$ converging locally uniformly to a function $u$ satisfying \eqref{eq unif holder}. To alleviate the notation, this subsequence is still denoted $(u^n)_{n\ge1}$. Recalling that $Y^n_t = u^n(t,X_t)$ and introducing $Y_t := u(t,X_t)$, for $t \in [0,T]$, we observe that
\begin{align}\label{eq conv Y}
\esp{\sup_{t \in [0,T]}|Y^n_t - Y_t|^2} \underset{n \rightarrow +\infty}{\longrightarrow} 0 \,,
\end{align}
since $t \mapsto (t,X_t)$ is a.s. continuous and $\{|Y^n|\}$ is bounded uniformly by a constant, see Lemma \ref{pr penal pseudo distance}.


\medskip

We conclude this section with the following lemma, which is used in the next section. This lemma provides a uniform upper bound on the second moment of the auxiliary process 
\begin{align*}
&\Gamma^{n,m}_t := 
\exp\left(  \int_0^t \left(1 + |\dot K^n_s| + |\dot K^m_s| \right) \ud s \right),\quad t\in[0,T],
\end{align*}
where we recall \eqref{eq.Kn.def}.

\begin{Lemma}
\label{cor:eps:sliceable}
Under Assumptions \ref{ass:intro:SDE} and \ref{ass:markov}, for any $\varepsilon>0$, there exists $N\geq1$ (independent of $n$) such that, for all $n\geq1$ and all $0 \leqslant k < N$, we have a.s.:
$$
\mathbb{E}_{t} \left[\int_t^{T(k+1)/N} |Z_s^{n}|^2+|\dot{K}^n_s| \ud s\right]\leqslant \varepsilon, \quad \forall t \in [Tk/N,T(k+1)/N].
$$
In particular, for any $\beta>0$, there exists a constant $C=C(\beta)$ (independent of $(n,m)$), s.t. 
$$
\mathbb{E}[(\Gamma^{n,m}_T)^{\textcolor{black}{\beta}}] \leqslant C,
$$
for all $n,m\geq1$.
\end{Lemma}

\begin{Remark}\label{rem:C.indep.x.1}
It is worth noticing that the constant $C$, appearing in Lemma \ref{cor:eps:sliceable}, does not depend on the initial value $x$ of the diffusion $X$, as follows from the proof of the lemma.
\end{Remark}

\proof
The proof of the first statement of the lemma is an improvement of the estimates in the proof of Lemma \ref{pr first estimates}, with the use of Corollary \ref{cor:unifHolder}.
We fix $t<t' \in [0, T]$, $\beta'>0$ and $\alpha>0$, and apply It\^o's formula to the process $(e^{\beta'(s-t)}(\alpha|Y^n_s|^2+(\phi_{\mathcal{C}}(Y^n_s)-\phi_{\mathcal{C}}(0))^2)_{s\in[t,t']}$ (recall \eqref{eq de penalised equation}) to obtain, as in the proof of Lemma \ref{pr first estimates},
\begin{align*}
 |Y^n_t|^2 + c \mathbb{E}_t\left[\int_t^{t'} |Z_s^n|^2+|\dot{K}^n_s| \ud s \right] &\leqslant \mathbb{E}_{t}\left[e^{\beta'(t'-t)}|Y^n_{t'}|^2 + C\int_t^{t'} e^{\beta'(s-t)} (1+|F(s,X_s,0,0)|^2) \ud s\right],
\end{align*}
which holds for large enough $\beta'$ and small enough $\alpha$.

Then, by using the upper bounds on $|Y^n|$, see Proposition \ref{pr penal pseudo distance 2} and on $|F(.,.,0,0)|$, see Assumption \ref{ass:markov}, we obtain: 
\begin{align*}
&\mathbb{E}_t\left[\int_t^{t'} |Z_s^n|^2+|\dot{K}^n_s| \ud s \right]\\
\leqslant& \mathbb{E}_{t}\left[e^{\beta'(t'-t)}|Y^n_{t'}|^2-|Y^n_t|^2 + C\int_t^{t'} e^{\beta'(s-t)} (1+|F(s,X_s,0,0)|^2) \ud s\right]\\
\leqslant& \mathbb{E}_{t}\left[(e^{\beta'(t'-t)}-1)|Y^n_{t'}|^2+ |Y^n_{t'}+Y^n_t||u^n(t',X_{t'})-u^n(t,X_t)|+\frac{C}{\beta'}(e^{\beta'(t'-t)} -1)\right]\\
\leqslant& C(\beta')(t'-t)+ C\mathbb{E}_{t}\left[(t'-t)^{\alpha'/2}+|X_{t'}-X_t|^{\alpha'}\right]
\leqslant C'(\beta')(t'-t)^{\alpha'/2},
\end{align*}
where $C'$ is independent of $n$, and we made use of Jensen's inequality and of standard SDE estimates on $X$ in the last inequality.
The above proves the first statement of the lemma.

\smallskip

To prove the second statement, we fix an arbitrary $\beta>0$ and consider $N$ corresponding to $\varepsilon=1/(8\beta)$. Then, the first statement of the lemma and the John-Nirenberg inequality yield:
\begin{align*}
\mathbb{E}\left[e^{2\beta \int_0^T |\dot{K}^n_s| + |\dot{K}^m_s| \ud s}\right] 
&\leqslant  \mathbb{E}\left[e^{2\beta \int_0^{T(N-1)/N} |\dot{K}^n_s| + |\dot{K}^m_s| \ud s}\mathbb{E}_{T(N-1)/N}\left[e^{2\beta \int_{T(N-1)/N}^T |\dot{K}^n_s|+ |\dot{K}^m_s| \ud s}\right] \right]\\
&\leqslant 2\mathbb{E}\left[e^{2\beta \int_0^{T(N-1)/N} |\dot{K}^n_s|+ |\dot{K}^m_s| \ud s} \right].
\end{align*}
Iterating the above, we obtain the desired estimate.
\eproof



\subsection{Existence and uniqueness}

We denote by $\{(Y^n,Z^n)\}_{n\geq1}$ a sequence of Markovian solutions to \eqref{eq de penalised equation} satisfying \eqref{eq conv Y} (whose existence is established in the previous subsection).
The goal of this subsection is to establish that $\{(Y^n,Z^n,K^n)\in\SP{2}\times \HP{2} \times \mathscr{K}^1\}_{n\geq1}$ (with $K^n$ defined in \eqref{eq.Kn.def})\footnote{The fact that $K^n\in\mathscr{K}^1$ follows from the inequality \eqref{eq.aprioriEst1.eq3} and the second statement of Proposition \ref{prop:penal.WellPosed}.} converges to a solution of the reflected BSDE \eqref{eq reflected bsde} and that this solution is unique in the appropriate class.

\begin{Theorem}
\label{thm:existence:Markov}
Let Assumptions \ref{ass:intro:SDE} and \ref{ass:markov} hold. Then, there exists a triplet $(Y,Z,K)\in\SP{2}\times \HP{2} \times \mathscr{K}^1$, such that
$$
\lim_{n\rightarrow\infty}(\|Y^n-Y\|_{\SP{2}},\|Z^n-Z\|_{ \HP{2}},\|K^n-K\|_{\SP{2}})=0\;.
$$
which solves \eqref{eq reflected bsde}. The process $K$ is absolutely continuous and satisfies, for all $\beta > 0$,
\begin{align}\label{eq exp moment}
\esp{e^{\beta\mathrm{Var}_T(K)}} < \infty\;.
\end{align}
Moreover, this solution is unique in the class $\mathscr{U}(1)$ (recall Definition \ref{de uniqueness class}).
\end{Theorem}

\begin{Remark}
If, in addition to Assumptions \ref{ass:intro:SDE} and \ref{ass:markov}, $g$ and $F$ are globally Lipschitz in $x$ (i.e., $\alpha=1$ in Assumption \ref{ass:markov}), then there exists a constant $C$ such that 
\begin{equation*}
|Z| \leqslant C, \quad \ud t \times \ud \mathbb{P}\text{-a.e.}
\end{equation*}
Indeed, using the same arguments as in the proof of Corollary \ref{cor:unifHolder}, we conclude that the conditions of Theorem 2.16 in \cite{Harter-Richou-19} are satisfied. The latter theorem yields the existence of a constant $C$, s.t. $|Z^n_t|\leq C$ for a.e. $(t,\omega)$ and for all $n$. Then, it follows that $|Z|\leq C$.
\end{Remark}

\begin{Remark}\label{rem:C.indep.x.2}
It is worth noticing that every exponential moment of $\mathrm{Var}_T(K)$ can be bounded by a constant that does not depend on the initial value $x$ of the diffusion $X$, as follows from Remark \ref{rem:C.indep.x.1} and from the proof of Theorem \ref{thm:existence:Markov}.
\end{Remark}


\proof 
1.a We first prove the uniqueness of the solution in the desired class. For any solution $(Y',Z',K')$ in $\mathscr{U}(1)$, we have
\begin{equation}\label{eq.Markov.weelPosed.Pf.eq1}
\mathbb{E}\left[ e^{\frac{p'}{R_0} \mathrm{Var}_T(K')} \right] <+\infty,
\end{equation}
for some $p'>1$.
Setting $1<p := (1+p')/2<p'$, $q' := p'/p>1$ and $q=q'/(q'-1)$, we obtain, using H\"older inequality,
\begin{align}
\esp{e^{\frac{p}{R_0} \left(\mathrm{Var}_T(K) + \mathrm{Var}_T(K') \right)}}
\le 
\esp{e^{\frac{q p}{R_0} \mathrm{Var}_T(K)} }^{\frac1{q}} \esp{e^{\frac{q' p}{R_0} \mathrm{Var}_T(K')} }^{\frac1{q'}}.
\end{align}
By \eqref{eq.Markov.weelPosed.Pf.eq1}, we have $\esp{e^{\frac{q' p}{R_0} \mathrm{Var}_T(K)} } = \esp{e^{\frac{p'}{R_0} \mathrm{Var}_T(K)} }<\infty$. Then using \eqref{eq exp moment}, which is proved below, we obtain
 \begin{align}
\esp{e^{\frac{p}{R_0} \left(\mathrm{Var}_T(K) + \mathrm{Var}_T(K') \right)}} <+\infty.
\end{align}
Proposition \ref{prop:stability2}, then, yields the uniqueness stated in the theorem.

\smallskip

\noindent 1.b The fact that $K$ is absolutely continuous is proved in Lemma \ref{le abs cont K}.

\smallskip

\noindent 2. Turning to the  existence part of the proof, we have already obtained the convergence of $\{Y^n\}$ -- recall \eqref{eq conv Y}. Moreover, it follows easily from Proposition \ref{pr penal pseudo distance 2} that, with probability one, $Y_t$ takes values in $\bar\cD$ for all $t\in[0,T]$.

We now turn to the convergence of $\{Z^n\}$. For $n,m \ge 1$, we denote 
$$
\df_t := f(t,Y^n_t,Z^n_t) - f(t,Y^m_t,Z^m_t),\quad \delta K :=  K^n- K^m,
$$ 
$$
\dY := Y^n - Y^m,\quad \dZ = Z^n - Z^m.
$$
Applying It\^o's formula to $(e^{\beta' s}|\dY_s|^2)_{s\in[t,T]}$, we obtain
\begin{align}
&|\dY_t|^2 + \int_t^T\!\!  e^{\beta' (s-t)} |\dZ_s|^2 \ud s = 2 \int_t^T\!\!  e^{\beta' (s-t)} \dY_s\cdot \df_s \ud s
- 2 \int_t^T\!\!  e^{\beta' (s-t)} \dY_s \cdot \delta \dot{K}_s \ud s\nonumber\\
&- 2 \int_t^T \!\! e^{\beta' (s-t)} \dY_s\cdot \dZ_s \ud W_s
- \beta' \int_t^T e^{\beta' (s-t)}|\dY_s|^2 \ud s.
\label{eq ito for dummies}
\end{align}
Choosing a large enough $\beta'>0$ and using the standard estimates, we deduce
\begin{align}\label{eq ref interm step}
\esp{\int_0^T\!\!  |\dZ_s|^2 \ud s} \le C \esp{\int_0^T\!\!  \left|\dY_s \cdot \delta \dot{K}_s \right|\ud s}.
\end{align}
Note that Lemma \ref{cor:eps:sliceable} yields the existence of a constant $C$, s.t.
$\esp{ \left(\int_0^T |\delta \dot{K}_s| \ud s\right)^2} \le C$, for all $n,m$.
Then, using the Cauchy-Schwartz inequality, we obtain
\begin{align*}
\esp{ \int_0^T\!\!  |\dY_s \cdot \delta \dot{K}_s| \ud s}
\le \esp{ \sup_{s\in[0,T]} |\dY_s|^2 }^\frac12 \left[\EE\left(\int_0^T|\delta \dot{K}_s| \ud s\right)^2\right]^\frac12.
\end{align*}
The above estimate, along with \eqref{eq ref interm step} and \eqref{eq conv Y}, implies that $\{Z^n\}_{n\ge 1}$ is a Cauchy sequence. 
Thus, there exists $(Y,Z)\in\SP{2}\times \HP{2}$ such that $(Y^n,Z^n) \rightarrow (Y,Z)$.

\smallskip

Next, we recall that
\begin{align*}
K^n_t = Y^n_t - Y^n_0 + \int_0^t f(s,Y^n_s,Z^n_s) \ud s - M^n_t,\quad M^n_t:= \int_0^t Z^n_s \ud W_s.
\end{align*}
Doob's maximal inequality implies that $\{M^n\}$ converges in $\SP{2}$ to $M$, with $M_t:= \int_0^t Z_s \ud W_s$.
As $f(t,\cdot,\cdot)$ is Lipschitz, we conclude that
\begin{equation}\label{eq.Kn.conv.K}
\|K^n-K\|_{\SP{2}}\rightarrow0,
\end{equation}
with the continuous process $K$ defined as
\begin{align*}
K_t := Y_t - Y_0 + \int_0^t f(s,Y_s,Z_s) \ud s - \int_0^t Z_s \ud W_s.
\end{align*}

\medskip

We now prove that $K \in \mathscr{K}^1$, and that $\ud K_t$ is directed along $\mathfrak{n}$ and is active only when $Y$ touches the boundary.
To this end, we define the auxiliary nondecreasing processes
$$
\hat K^n_t:=\int_0^t n\psi(Y^n_s)(1+|Z^n_s|^2) \ud s,\quad t\in[0,T].
$$
From Lemma \ref{pr first estimates} we deduce the existence of a constant $C$, s.t. $\esp{\hat K^n_T} \le C$ for all $n$.
Then, using Proposition 3.4 in \cite{campi2006super}, we know that there exists a nondecreasing nonnegative process $\hat K$, two sequences of integers $\{p\leq N_p\}$, with $p\rightarrow\infty$, and a family of numbers $\{\lambda^p_r\}$, with $\sum_{r=p}^{N_p}\lambda^p_r =1$, such that
\begin{align}
\P\left( {}^p\!{\hat K}_t := \sum_{r=p}^{N_p}\lambda^p_r \hat K^r_t \rightarrow \hat K_t\,,\,\forall t \in[0,T]\right) = 1\,.
\end{align}
The above implies that the measure induced by $\ud\, {}^p\!{\hat K}_t$ on $[0,T]$ converges a.s. to $\ud \hat K_t$.
Then, for any bounded continuous process $\chi$ and any $0\leq t_1<t_2\leq T$,
\begin{align}
\eta^{p}(t_1,t_2)
:=\int_{t_1}^{t_2} \chi_t \sum_{r=p}^{N_p}\lambda^p_r r\psi(Y^r_t)(1+|Z^r_t|^2) \ud t
=\int_{t_1}^{t_2} \chi_t \ud\, {}^p\!{\hat K}_t \rightarrow \int_0^T \chi_t \ud \hat K_t,
\quad a.s.\label{eq.existence.etap.conv}
\end{align}
From the first statement of Lemma \ref{cor:eps:sliceable} (with the use of Proposition \ref{pr penal pseudo distance 2}), we conclude that, for any $\varepsilon>0$, there exists $N\geq1$ (independent of $n$) such that, for all $p$ and all $0 \leqslant k < N$, we have a.s.:
$$
\EFp{t}{ |\eta^p(t,T(k+1)/N)|}\leqslant \varepsilon, \quad \forall t \in [Tk/N,T(k+1)/N].
$$
Then, repeating the last part of the proof of Lemma \ref{cor:eps:sliceable}, we conclude that, for any $\beta>0$, there exists a constant $C$, s.t.
$$
\esp{e^{\beta\eta^p(0,T)} }\leq C,\quad \forall\,p.
$$
Thus, the family $\{\exp(\beta\eta^p(0,T))\}_p$ is uniformly integrable. The latter implies, in particular, that the convergence in \eqref{eq.existence.etap.conv} holds in $\mathcal{L}^1$ and that all exponential moments of $\hat K_T$ are finite.

\smallskip

Next, we define
$$
{}^p\!{K}_t := \sum_{r=p}^{N_p}\lambda^p_r K^r_t,\quad t\in[0,T].
$$
We also denote by $\widetilde{\nabla\psi}$ a Lipschitz extension of $\nabla\psi$ into $\cD$ (constructed as in Remark \ref{rem:explain.Ito}).
Then, for any event $A$ and any $t\in[0,T]$, we have:
$$
\esp{ {}^p\!{K}_t\, \bone_A }= \esp{\int_0^t \sum_{r=p}^{N_p}\lambda^p_r \widetilde{\nabla\psi}(Y^r_s) r\psi(Y^r_s)(1+|Z^r_s|^2) \ud s\,\bone_A}
$$
$$
= \esp{ \int_0^t \widetilde{\nabla\psi}(Y_s) \sum_{r=p}^{N_p}\lambda^p_r r\psi(Y^r_s)(1+|Z^r_s|^2) \ud s\,\bone_A}
$$
$$
+ O\left(\esp{ \sum_{r=p}^{N_p}\lambda^p_r \sup_{s\in[0,t]}|Y^r_s-Y_s| \int_0^t (1+|Z^r_s|^2) \ud s} \right)
$$
$$
\rightarrow \esp{ \int_0^t \widetilde{\nabla\psi}(Y_s) \ud \hat K_s \ud s\,\bone_A},
$$
where we used \eqref{eq.existence.etap.conv}, and its $\mathcal{L}^1$ version, and the estimate
$$
\esp{ \sup_{s\in[0,t]}|Y^r_s-Y_s| \int_0^t (1+|Z^r_s|^2) \ud s}
\leq C \|Y^r-Y\|_{\SP{2}} \esp{ e^{\int_0^t (1+|Z^r_s|^2)} }^{1/2}
\leq C \|Y^r-Y\|_{\SP{2}},
$$
which follows from Lemma \ref{cor:eps:sliceable}.

On the other hand, as $K^n$ converges to $K$ in $\SP{2}$, $\esp{ {}^p\!{K}_t\, \bone_A}$ converges to $\esp{ K_t\bone_A}$,
and since $A$ is arbitrary and $K_\cdot$ is continuous, we conclude:
\begin{equation}\label{eq.pushOrt}
\P\left(K_t = \int_0^t \widetilde{\nabla\psi}(Y_s) \ud \hat K_s \ud s,\quad \forall\,t\in[0,T]\right)=1.
\end{equation}
Note that the integrability of $\hat K_T$ and the above representation, in particular, imply $K \in \mathscr{K}^1$.

\smallskip

It only remains to show that
\begin{equation}\label{eq.pushAtBndry}
\int_0^T \bone_{\cD}(Y_t) \ud \hat K_t=0.
\end{equation}
To this end, we choose an arbitrary Lipschitz $f$ supported in $\cD$ and any event $A$, to obtain:
$$
\esp{ \int_0^T f(Y_t) \ud \hat K_t\,\bone_A }
= \lim_{n\rightarrow\infty} \esp{ \int_0^T f(Y_t) \sum_{r=p}^{N_p}\lambda^p_r r\psi(Y^r_t)(1+|Z^r_t|^2) \ud t\,\bone_A}
$$
$$
= \lim_{n\rightarrow\infty} \esp{ \int_0^T \sum_{r=p}^{N_p}\lambda^p_r f(Y^r_t) r\psi(Y^r_t)(1+|Z^r_t|^2) \ud t\,\bone_A}
$$
$$
+ O\left(\esp{\sum_{r=p}^{N_p}\lambda^p_r \sup_{t\in[0,T]}|Y^r_t-Y_t| \int_0^T (1+|Z^r_t|^2) \ud t}\right)
=0.
$$
As $A$ is arbitrary, we conclude that, for any Lipschitz $f$ supported in $\cD$, we have $\int_0^T f(Y_t) \ud \hat K_t=0$ a.s..
Approximating $\bone_{\cD}$ with a sequence of such $f$, we, e.g., use the monotone convergence theorem to deduce \eqref{eq.pushAtBndry}.
Combining the latter with \eqref{eq.pushOrt}, we obtain \eqref{eq reflected bsde}(ii) and conclude the proof of the first part of Theorem \ref{thm:existence:Markov}. 
\eproof

\smallskip

\begin{Remark}
\label{rem:unique.Markov}
Theorem \ref{thm:existence:Markov} implies that, under Assumptions \ref{ass:intro:SDE}, \ref{ass:markov} and \ref{ass:general} with $\theta = 1$, there exists a unique solution to \eqref{eq reflected bsde} in $(Y,Z,K)\in\SP{2}\times \HP{2} \times \mathscr{K}^1$.
\end{Remark}



\section{Well-posedness beyond Markovian framework}
\label{se:nonMarkov}

\subsection{Discrete path-dependent framework}

In this subsection, we extend the existence and uniqueness result obtained in a Markovian framework, see Theorem \ref{thm:existence:Markov}, to a discrete path-dependent framework.

\begin{Assumption}
 \label{ass:pathdependent}
Let $\ell$ be an arbitrary strictly positive integer and consider the partition $0 =t_0 < t_1<...< t_\ell=T$ of $[0,T]$. We assume that
$$
 \xi =g(X_{t_1},...,X_{t_\ell}) \quad \text{and} \quad f(s,y,z)=F(s,X_{t_1 \wedge s},\cdots,X_{t_\ell \wedge s},y,z)$$
where 
\begin{enumerate}[(i)]
 \item $g$ is $\alpha$-H\"older and takes values in $\bar{\cD}$,
 \item $F$ is measurable in all variables, globally Lipschitz in $(y,z)$, uniformly over $(x_1,...,x_{\ell})$, and globally $\alpha$-H\"older in $(x_1,...,x_{\ell})$, uniformly over $(y,z)$, and $|F(\cdot,\cdots,0,0)|$ is bounded.
\end{enumerate}
\end{Assumption}

We note that $\ell=1$ corresponds to the Markovian framework with an extra regularity assumption on the generator with respect to $x$.
We also recall that Assumptions \ref{ass:intro.domain} and \ref{ass:intro.data} hold throughout Section \ref{se:nonMarkov} even if not cited explicitly.

\begin{Theorem}
 \label{thm:pathdependent}
 Let Assumptions \ref{ass:intro:SDE}  and Assumption \ref{ass:pathdependent} hold. Then, there exists a triplet $(Y,Z,K)\in\SP{2}\times \HP{2} \times \mathscr{K}^1$ that solves \eqref{eq reflected bsde}. Moreover, all exponential moments of $\mathrm{Var}_T(K)$ are finite, and this solution is unique in the class $\mathscr{U}(1)$ (recall Definition \ref{de uniqueness class}).
\end{Theorem}
\proof
Once the finiteness of the exponential moments of $\mathrm{Var}_T(K)$ is proven, the uniqueness of the solution in the class $\mathscr{U}(1)$ follows from the same arguments as in step 1.a of the proof of Theorem \ref{thm:existence:Markov}.
Let us now prove the existence part of the theorem. To this end, we use the backward recursion to construct a solution on each interval $[t_i,t_{i+1}]$ for $0 \leqslant i \leqslant \ell -1$.

Since the case $\ell =1$ corresponds to the already treated Markovian framework, we assume that $\ell >1$ and we consider the interval time $[t_{\ell -1},T]$.
For any $(t,x) \in [0,T] \times \R^{d'}$, we denote by $X^{t,x}$ the unique solution of \eqref{SDE} on $[t,T]$, which starts from $x$ at time $t$. We reserve the notation $X$ for the original diffusion started at time zero. For any $\mathbf{x}=(\mathbf{x}_1,...,\mathbf{x}_{\ell-1}) \in (\R^{d'})^{\ell-1}$, we denote by $(Y^{\mathbf{x}},Z^{\mathbf{x}},K^{\mathbf{x}})$ the solution of \eqref{eq reflected bsde} on $[t_{\ell-1},T]$, with the terminal condition $g(\mathbf{x},X_T^{t_{\ell-1},\mathbf{x}_{\ell-1}})$ and with the generator $F(.,\mathbf{x},X_.^{t_{\ell-1},\mathbf{x}_{\ell-1}},.,.)$, whose existence is ensured by Theorem \ref{thm:existence:Markov} and whose uniqueness in an appropriate class follows from Theorem \ref{th:uniqueness}. 

Next, we denote by $(Y^{\mathbf{x},n},Z^{\mathbf{x},n})$ a Markovian solution of the penalized BSDE \eqref{eq de penalised equation} on $[t_{\ell-1},T]$, whose existence is ensured by Proposition \ref{prop:penal.WellPosed}. In particular, there exist measurable functions $u^n(\mathbf{x},.,.)$ and $v^n(\mathbf{x},.,.)$ such that  
\begin{equation*}
Y^{\mathbf{x},n}_t = u^n(\mathbf{x},t,X_t^{t_{\ell-1},\mathbf{x}_{\ell-1}}), \quad  Z^{\mathbf{x},n}_t = v^n(\mathbf{x},t,X_t^{t_{\ell-1},\mathbf{x}_{\ell-1}}).
\end{equation*}
By considering a sequence of Lipschitz approximations of \eqref{eq de penalised equation}, given by \eqref{eq.penalized.trunc}, we apply Theorem 5.4 in \cite{Hu-Ma-04} and, passing to the limit for the Lipschitz approximations as in the proof of Proposition \ref{prop:penal.WellPosed}, we conclude that a Markovian solution to \eqref{eq de penalised equation} can be constructed so that $u^n$ and $v^n$ are jointly measurable in all variables.
Passing to the limit in $n$ along a subsequence, we use Theorem \ref{thm:existence:Markov} and the uniform H\"older estimate in Corollary \ref{cor:unifHolder}, to deduce the existence of jointly measurable functions $u$ and $v$ satisfying 
\begin{equation}\label{eq.nonMarkov.disc.thm1.eq1}
Y^{\mathbf{x}}_t = u(\mathbf{x},t,X_t^{t_{\ell-1},\mathbf{x}_{\ell-1}}), \quad  Z^{\mathbf{x}}_t = v(\mathbf{x},t,X_t^{t_{\ell-1},\mathbf{x}_{\ell-1}}).
\end{equation}

Then, by denoting $\mathbf{X}=(X_{t_1},...X_{t_{\ell-1}})$, we consider the progressively measurable processes $(Y^{\mathbf{X}}_t,Z^{\mathbf{X}}_t)_{t \in [t_{\ell-1},T]}$ and define
$$K_t^{\mathbf{X}}:= Y_t^{\mathbf{X}}-Y_{t_{\ell-1}}^{\mathbf{X}} +\int_{t_{\ell-1}}^t F(s,\mathbf{X},X_s^{t_{\ell-1},X_{t_{\ell-1}}},Y_s^{\mathbf{X}},Z_s^{\mathbf{X}})\ud s - \int_{t_{\ell-1}}^t Z_s^{\mathbf{X}}\ud W_s,\quad t_{\ell-1} \leqslant t \leqslant T.$$
We note that $X_s^{t_{\ell-1},X_{t_{\ell-1}}}=X_s$ and that $(Y^{\mathbf{X}}_t,Z^{\mathbf{X}}_t,K^{\mathbf{X}}_t)_{t \in [t_{\ell-1},T]}$ is a solution of \eqref{eq reflected bsde} on the time interval $[t_{\ell-1},T]$ satisfying $K_{t_{\ell-1}}^{\mathbf{X}}=0$. 

\smallskip

In order to iterate this construction and to extend the solution to the time interval $[t_{\ell-2},t_{\ell-1}]$, we have to ensure that the associated terminal condition $Y^{\mathbf{X}}_{t_{\ell-1}}$ of the reflected BSDE \eqref{eq reflected bsde} on $[t_{\ell-2},t_{\ell-1}]$ is an $\alpha$-H\"older function of $\mathbf{X}$. To this end, we recall the function $u$ in \eqref{eq.nonMarkov.disc.thm1.eq1} and define, for all $\tilde{\mathbf{x}}= (\mathbf{x}_1,...,\mathbf{x}_{\ell-2}) \in (\R^{d'})^{\ell-2}$ and $\mathbf{x}_{\ell-1} \in \R^{d'}$, the deterministic function
$$
\tilde{g}(\tilde{\mathbf{x}},\mathbf{x}_{\ell-1}):=u(\tilde{\mathbf{x}},\mathbf{x}_{\ell-1},t_{\ell-1},\mathbf{x}_{\ell-1})
= Y^{\tilde{\mathbf{x}},\mathbf{x}_{\ell-1}}_{t_{\ell-1}}.
$$
Let us prove that this function is $\alpha$-H\"older. Indeed, for any $\mathbf{x}:=(\tilde{\mathbf{x}},\mathbf{x}_{\ell-1})\in(\R^{d'})^{\ell-1}$ and $\mathbf{x}':=(\tilde{\mathbf{x}}',\mathbf{x}_{\ell-1}')\in(\R^{d'})^{\ell-1}$, Proposition \ref{prop:stability2} with $p=2$ yields
\begin{align*}
 &|\tilde{g}(\mathbf{x})-\tilde{g}(\mathbf{x}')| \leqslant \| Y^{\mathbf{x}} - Y^{\mathbf{x}'} \|_{\SP{2}}\\
 \leqslant& C\mathbb{E}\left[ \left|g(\mathbf{x},X_T^{t_{\ell-1},\mathbf{x}_{\ell-1}})- g(\mathbf{x}',X_T^{t_{\ell-1},\mathbf{x}_{\ell-1}'})\right|^4\right]^{1/4}\\
 &+C \mathbb{E}\left[\left(\int_{0}^T |F(s,\mathbf{x},X_s^{t_{\ell-1},\mathbf{x}_{\ell-1}}, Y_s^{\mathbf{x}},Z_s^{\mathbf{x}})-F(s,\mathbf{x}',X_s^{t_{\ell-1},\mathbf{x}_{\ell-1}'}, Y_s^{\mathbf{x}},Z_s^{\mathbf{x}})| \ud s\right)^4\right]^{1/4}\\
\leqslant & C\left( |\mathbf{x}-\mathbf{x}'|^{\alpha}+ \mathbb{E}\left[\sup_{0 \leqslant s \leqslant T} |X_s^{t_{\ell-1},\mathbf{x}_{\ell-1}}-X_s^{t_{\ell-1},\mathbf{x}_{\ell-1}'}|^{4\alpha}\right]^{1/4}\right),
\end{align*}
with a constant $C$ that does not depend on $x$ (see Remark \ref{rem:C.indep.x.2}).
Then, the Jensen's inequality and the standard SDE estimates yield 
$$
|\tilde{g}(\mathbf{x})-\tilde{g}(\mathbf{x}')| \leqslant  C |\mathbf{x}-\mathbf{x}'|^{\alpha},
$$
which gives us the $\alpha$-H\"older property of $\tilde g$. Considering the reflected BSDE \eqref{eq reflected bsde} on $[t_{\ell-2},t_{\ell-1}]$, with the terminal condition $Y^{\mathbf{X}}_{t_{\ell-1}}=\tilde{g}(X_{t_1},...,X_{t_{\ell-1}})$ and with the generator 
$$
F(s,X_{t_1 \wedge s},...,X_{t_{\ell-2} \wedge s},X_{t_{\ell-1} \wedge s},X_{t_{\ell-1} \wedge s},y,z),
$$
we deduce, as in the first part of the proof, that it has a solution in the form \eqref{eq.nonMarkov.disc.thm1.eq1}.

\smallskip

Finally, iterating the above construction, we concatenate the ``$Y$'' and ``$Z$'' parts of the solutions constructed in individual sub-intervals, and we sum up the ``$K$'' parts (assuming that every individual ``$K$'' part is extended continuously as a constant to the left and to the right of the associated sub-interval).
It is easy to see that the resulting process $(Y,Z,\tilde{K}) \in \SP{2} \times \HP{2} \times \mathscr{K}^1$ is a solution of \eqref{eq reflected bsde} on $[0,T]$.
\eproof

\subsection{General case}


\begin{Theorem}
 \label{thm:general}
 Let Assumption \ref{ass:general} hold with $\theta = 2$. Then, there exists a triplet $(Y,Z,K)\in\SP{2}\times \HP{2} \times \mathscr{K}^1$ that solves \eqref{eq reflected bsde}, and this solution is unique in the class $\mathscr{U}(2)$.
\end{Theorem}

\proof
The uniqueness part of the theorem is a direct consequence of Proposition \ref{A priori estimate reflected BSDE} and Corollary \ref{th:uniqueness}.
Let us now prove the existence part: To do so, we shall construct a Cauchy sequence of approximating reflected BSDEs.
\\
First, we observe that the terminal condition $\xi$ can be approximated by a sequence of random
variables of the form $\xi^n:=g_n(W_{t_1},...,W_{t_n})$, where $g_n$ is infinitely differentiable. The sequence $(\xi^n)_{n \in \mathbb{N}^*}$ can be chosen so that it converges to $\xi$ in $\mathscr{L}^q$, for any $q \geqslant 1$ (see, e.g., \cite{Nualart-06}).
In particular,
\begin{equation}\label{eq.SecGen.mainThm.xin.approx.xi}
\lim_{n\rightarrow\infty}\EE \left[|\xi-\xi^n|^{2p/(p-1)}\right] = 0,
\end{equation}
with $p>1$ appearing in Proposition \ref{A priori estimate reflected BSDE}.
Replacing $g_n$ by $g_n \wedge \|\xi\|_{\mathscr{L}^{\infty}}$, we can assume $\|\xi^n\|_{\mathscr{L}^{\infty}} \leqslant \|\xi\|_{\mathscr{L}^{\infty}}$. We observe that $\xi^n$ satisfies Assumption \ref{ass:pathdependent}(i) for $X=W$.

Second, to approximate the generator, for every $n \in \mathbb{N}^*$, we denote by $\mathcal{K}^n$ the closed ball in $\R^{d\times d'}$ of radius $n$ centered at zero, and choose a sequence of numbers $\epsilon_n\downarrow0$. 
We set 
\begin{align*}
\ell_n := \| f(\cdot,0,0)\|_{\mathscr{L}^{\infty}} + K_{f,y}\sup_{y \in \bar{\cD}}|y| + nK_{f,z},
\end{align*}
recalling Assumption \ref{ass:intro.data}.
For each $n$, we introduce  $\mathfrak{L}^{n}$ the compact convex subset of $\mathscr{C}(\bar\cD\times\mathcal{K}^n)$ (the space of continuous function endowed with the uniform norm denoted $\|\cdot\|_s$) consisting of Lipschitz functions, with the Lipschitz coefficients $K_{f,y}$ in the variable $y \in \bar\cD$ and $K_{f,z}$ in the variable $z\in \mathcal{K}^n$, and whose (uniform) norm is bounded by $\ell_n$. 
Note that $f_{|\mathcal{K}^n}$ is valued in $\mathfrak{L}^n$. We are now going to build an approximation of $f_{|\bar\cD\times\mathcal{K}^n}$
in $\mathfrak{L}^n$ satisfying Assumption \ref{ass:pathdependent} (for $X=W$). \\
Let  $\{\phi_{n}^{m}\}_{m=1}^{M_n}$ be an $\epsilon_n$-cover of the compact set $\mathfrak{L}^n$ with $M_n$ a positive integer. 
We denote by $\tilde f^{n}(t,\cdot)$ the (measurable selection of the) proximal projection of $f_{|\bar{\cD}\times\cK^n}(t,\cdot)$ on $\{\phi_{n}^{m}\}_{m=1}^{M_n}$ . It satisfies
\begin{align*}
\tilde f^{n}(t,\cdot) = \sum_{m=1}^{M_n}\phi^m_n(\cdot)(\tilde{\eta}^n_t)^m =: \phi_n(\cdot) \tilde{\eta}^n_t \quad \text{ and } \quad \|\tilde{f}^{n}(t,\cdot) - f_{|\bar{\cD}\times\cK^n}(t,\cdot)\|_s \le \epsilon_n \quad \text{a.s.}
\end{align*}
with $\tilde{\eta}^n$ a progressively measurable process which takes its value on the (non zero) extremal points of 
$\cS_{M_n} := \set{x \in \R^{M_n}\,|\, 0 \le x^m \le 1, \sum_{m=1}^{M_n} x^m \le 1}$. 
Then, using the dominated convergence theorem, we have
\begin{align}\label{eq step 1 approx}
\esp{ \int_0^T \|f_{|\bar{\cD}\times\cK^n}(t,\cdot) - \tilde f^{n}(t,\cdot)\|_s^{2p/(p-1)} \ud t }\le T(\epsilon_n)^{2p/(p-1)}.
\end{align}

We now classically approximate $(\tilde{\eta}^n_t)_{t\in[0,T]}$ by an adapted process 
$(\hat{\eta}^n_t)_{t\in[0,T]}$ piecewise constant on a grid $\Pi_n := \set{t_0 = 0 <\dots< t^n_k <\dots< t^n_{\kappa_n}=T}$.
This process can be chosen to be $\cS_{M_n}$-valued and satisfying
\begin{align*}
\esp{\int_0^T |\tilde{\eta}^n_t - \hat{\eta}^n_t|^{2p/(p-1)}\ud t} \le \frac{\epsilon_n}{(M_n\ell_n^2)^{p/(p-1)}} \;.
\end{align*}
Setting
\begin{align}
\hat{f}^n(t,\cdot) = \sum_{k=0}^{\kappa_n-1} \phi_n(\cdot) \hat\eta^{n}_{t^n_k} \bone_{(t^n_{k},t^n_{k+1}]}(t),
\end{align}
which is $\mathfrak{L}^n$-valued, as a random convex combination of  the $\set{\phi_n^m}_{m=1}^{M_n}$,
we compute
\begin{align}\label{eq step 2 approx}
\esp{ \int_0^T \|\tilde{f}^n(t,\cdot) -  \hat{f}^{n}(t,\cdot)\|_s^{2p/(p-1)} \ud t }\le \epsilon_n.
\end{align}

Then, we apply the approximation result of \cite{Nualart-06} for each $\hat\eta^{n}_{t^n_k}$.
Introducing, if necessary, a finer grid $\Re_n \subset \Pi_n$, we set
$$
\eta^{n}_{t^n_k}:=
\mathfrak{P}_{\cS} \left [ r^{n}_k\left( (W_{r})_{r \in \Re_n , r \le t^n_k}\right)\right],
$$
where  $r^n_k$ is a smooth function with values in $\R^{M_n}$ and $\mathfrak{P}_{\cS}$ the (orthogonal) projection onto $\cS_{M_n}$.
We can chose $r^n_k$ such that
\begin{align}
\esp{| \eta^{n}_{t_k} - \hat{\eta}^{n}_{t_k}|^{2p/(p-1)}} \le \frac{\epsilon_n}{(M_n\ell_n^2)^{p/(p-1)}} \;.
\end{align}
Setting $f^n(t,\cdot)=\sum_{k=0}^{\kappa_n-1} \phi_n(\cdot) \eta^{n}_{t^n_k} \bone_{(t^n_{k},t^n_{k+1}]}(t)$, which belongs to $\mathfrak{L}^n$, we have
\begin{align}\label{eq step 3 approx}
\esp{ \int_0^T \|\hat{f}^n(t,\cdot) -  {f}^{n}(t,\cdot)\|_s^{2p/(p-1)} \ud t }\le T\epsilon_n.
\end{align}

Collecting the above, we conclude that
\begin{equation}\label{eq.SecGen.mainThm.fn.approx.f}
\lim_{n\rightarrow\infty}\EE \int_0^T \sup_{y\in\bar\cD,\,z\in\mathcal{K}^n}|f(t,y,z) 
- f^n(t,y,z)|^{2p/(p-1)} \ud t =0.
\end{equation}
We extend $f^n(t,y,\cdot)$  to $\R^{d\times d'}\setminus\mathcal{K}^n$ as a constant in each radial direction, so that its uniform norm and Lipschitz coefficient do not change.

When, moreover, $f$ satisfies Assumption \ref{ass:general}-(i) (resp.  Assumption \ref{ass:general}-(iii)), the above construction allow to build an approximating sequence $f^n$ having the same properties. One simply
works with $\tilde{\mathfrak{L}}^n$ instead of $\mathfrak{L}^n$, where $\tilde{\mathfrak{L}}^n$  is the closed convex subset of $\mathfrak{L}^n$ whose function satisfies Assumption \ref{ass:general}-(i)  (resp.  Assumption \ref{ass:general}-(iii)).
\smallskip

Now, since, for any $n \in \mathbb{N}^*$, $\xi^n$ and $f^n$ satisfy Assumption \ref{ass:pathdependent}, we can invoke Theorem \ref{thm:pathdependent} to obtain an unique solution $(Y^n,Z^n,K^n)\in \SP{2} \times \HP{2} \times \mathscr{K}^1$ to \eqref{eq reflected bsde} associated with this data. Thanks to Proposition \ref{A priori estimate reflected BSDE}, we are allowed to apply Proposition \ref{prop:stability2}: for all $n,m \in \mathbb{N}^*$,
\begin{align}
 & \|Y^n-Y^m\|_{\SP{2}}+\|Z^n-Z^m\|_{\HP{2}}+\|K^n-K^m\|_{\SP{2}}\nonumber\\
 &\leqslant C\mathbb{E}\left[|\xi^n-\xi^m|^{2p/(p-1)}\right]^{(p-1)/(2p)}\label{eq.SecGen.mainThm.stability} \\
 &+ C\mathbb{E}\left[ \left(\int_0^T |f^n(s,Y_s^n,Z_s^n)-f^m(s,Y_s^n,Z_s^n)|\ud s\right)^{2p/(p-1)} \right]^{(p-1)/(2p)},\nonumber
\end{align}
with a constant $C$ that does not depend on $n$ and $m$. 

Applying the Cauchy-Schwartz, Jensen's, and Chebyshev's inequalities, we obtain
\begin{align*}
&\esp{\left(\int_0^T (1+|Z^n_t|)\1_{\set{|Z^n_t|>n}}\ud t\right)^{2p/(p-1)}} \\
&\leqslant T^{p/(p-1)-1/2}\esp{\left(\int_0^T (1+|Z^n_t|)^2 \ud t\right)^{2p/(p-1)}}^{1/2} \esp{ \int_0^T \1_{\set{|Z^n_t|>n}}\ud t}^{1/2} \\
&\leqslant \frac{T^{p/(p-1)-1/2}}{n}\esp{\left(\int_0^T (1+|Z^n_t|)^2 \ud t\right)^{2p/(p-1)}}^{1/2} \esp{\int_0^T |Z^n_t|^2 \ud t}^{1/2}.
\end{align*}
Using Proposition \ref{A priori estimate reflected BSDE} and the energy inequality for BMO martingales, we can estimate $\mathbb{E}\left(\int_0^T (1+|Z^n_t|)^2 \ud t\right)^{2p/(p-1)}$ uniformly over $n$.
Then, for all $m\geqslant n$, we obtain from the above estimate:
\begin{align}
&\mathbb{E}\left[\left(\int_0^T |f^n(t,Y^n_t,Z^n_t)-f^m(t,Y^n_t,Z^n_t)|\ud t\right)^{2p/(p-1)}\right]^{(p-1)/(2p)}\nonumber\\
&\leqslant C \mathbb{E}\left[\left(\int_0^T (1+|Z^n_t|)\1_{\set{|Z^n_t|>n}}\ud t\right)^{2p/(p-1)}\right]^{(p-1)/(2p)} \nonumber\\
&\quad + C\mathbb{E}\left[\left(\int_0^T |f^n(t,Y^n_t,Z^n_t)-f^m(t,Y^n_t,Z^n_t)|\1_{\set{|Z^n_t|\leqslant n}}\ud t\right)^{2p/(p-1)}\right]^{(p-1)/(2p)}
\label{eq.SecGen.mainThm.fn.fm.est}\\
&\leqslant \frac{C}{n^{(p-1)/(2p)}} +C\esp{\int_0^T |f^n(t,Y^n_t,Z^n_t)-f(t,Y^n_t,Z^n_t)|^{2p/(p-1)} \1_{\set{Y^n_t \in \bar{\cD},|Z^n_t|\leqslant n}}\ud t}^{(p-1)/(2p)}\nonumber\\
&\quad +C\esp{\int_0^T |f^m(t,Y^n_t,Z^n_t)-f(t,Y^n_t,Z^n_t)|^{2p/(p-1)} \1_{\set{Y^n_t \in \bar{\cD},|Z^n_t|\leqslant m}}\ud t}^{(p-1)/(2p)}.
\nonumber
\end{align}
In view of \eqref{eq.SecGen.mainThm.fn.approx.f}, the right hand side of the above vanishes as $n,m\rightarrow\infty$.
Collecting \eqref{eq.SecGen.mainThm.xin.approx.xi}, \eqref{eq.SecGen.mainThm.stability} and \eqref{eq.SecGen.mainThm.fn.fm.est}, we conclude:
\begin{align*}
 & \|Y^n-Y^m\|_{\SP{2}}+\|Z^n-Z^m\|_{\HP{2}}+\|K^n-K^m\|_{\SP{2}} \xrightarrow{n,m \rightarrow + \infty}0.
\end{align*}
In other words, $(Y^n,Z^n,K^n)_{n \in \mathbb{N}^*}$ is a Cauchy sequence in $\SP{2} \times \HP{2} \times \SP{2}$. Then, there exists $(Y,Z,K) \in \SP{2} \times \HP{2} \times \SP{2}$ such that $(Y^n,Z^n,K^n) \xrightarrow{n \rightarrow + \infty} (Y,Z,K)$. Moreover, $Y$ takes values on $\bar{\cD}$. As $(Y^n,Z^n,K^n)$ is the unique solution to \eqref{eq reflected bsde} associated to the terminal condition $\xi^n$ and the generator $f^n$, and since we have
\begin{align*}
&\mathbb{E}\left[\int_t^T |f^n(s,Y^n_s,Z_s^n)-f(s,Y_s,Z_s)|\ud s\right]\\
\leqslant & \mathbb{E}\left[\int_t^T |f^n(s,Y_s,Z_s)-f(s,Y_s,Z_s)|\ud s\right]+C\mathbb{E}\left[\int_t^T |Y^n_s-Y_s| + |Z^n_s-Z_s|\ud s\right],
\end{align*}
we can easily pass to the limit in \eqref{eq reflected bsde}(i) to show that $(Y,Z,K)$ satisfies \eqref{eq reflected bsde}(i).

\smallskip

It remains to prove that $K \in \mathscr{K}^1$, that $\ud K_t$ is directed along $\mathfrak{n}(Y_t)$ and that it is active only when $Y$ touches the boundary (the latter two properties will be shown via the alternative characterization given by Lemma \ref{lem alternative characterisation}).
Repeating the derivation of \eqref{ineq1}-\eqref{ineq2} for $(Y^n,Z^n,K^n)$, but without taking the conditional expectations, and with $\beta=0$, we obtain:
\begin{align*}
\int_0^{T} \ud \mathrm{Var}_s (K^n) \leqslant C \left(|\xi^n|^2+\int_0^{T}  2Y_s^n\cdot f(s,Y_s^n,Z_s^n) \ud s- \int_0^T 2Y_s^n Z_s^n \ud W_s\right),
\end{align*}
where the constant $C$ does not depend on $n$. The right hand side of the above inequality converges in probability, hence it also converges a.s. up to a subsequence which we still denote $\{(Y^n,Z^n,K^n)\}$. Then, $\{\textrm{Var}_T(K^n)\}_{n \in \mathbb{N}^*}$ is a.s. bounded uniformly over $n$, and Fatou's lemma yields that $\textrm{Var}_T(K)$ is a.s. bounded -- i.e. $K$ is a bounded variation process. Thanks to Proposition \ref{A priori estimate reflected BSDE}, $\{\textrm{Var}_T(K^n)\}_{n \in \mathbb{N}^*}$ is uniformly integrable and, hence, $K \in \mathscr{K}^1$. As $(Y^n,Z^n,K^n)$ solve \eqref{eq reflected bsde} with the terminal condition $\xi^n$ and the generator $f^n$, Lemma \ref{lem alternative characterisation} yields the existence of a constant $c$, independent of $n$, such that, for all continuous adapted process $V$ with values in $\bar{\cD}$, we have
$$\int_0^T (Y_s^n - V_s) \ud K_s^n + c|Y_s^n-V_s|^2 \mathfrak{n}(Y_s^n)\ud K_s^n \geqslant 0 \quad \textrm{a.s.}$$
Finally, we use Lemma 5.8 in \cite{Gegout-Petit-Pardoux-96} to pass to the limit in the above inequality, to obtain
$$\int_0^T (Y_s - V_s) \ud K_s + c|Y_s-V_s|^2 \mathfrak{n}(Y_s)\ud K_s \geqslant 0 \quad \textrm{a.s.}$$
and complete the proof of the theorem via another application of Lemma \ref{lem alternative characterisation}.
\eproof

\section{Connection to Brownian $\Gamma$-martingales}
\label{sect:sphere}

It turns out that the solutions to reflected BSDEs in non-convex domains, defined via \eqref{eq reflected bsde} and constructed in the previous sections, are naturally connected to the notion of martingales on manifolds (a.k.a. $\Gamma$-martingales -- see \cite{Emery-89}). In this section, we investigate this connection more closely, in particular, discovering a new proof of the existence and uniqueness of a Brownian martingale with a prescribed terminal value on a section of a sphere and showing the sharpness of the weak star-shape property in Assumption \ref{ass:intro.domain}.

\smallskip

The connection to martingales on manifolds is made precise by the following proposition, which states that, under certain assumptions, the $Y$-component of the solution to \eqref{eq reflected bsde} always stays on the boundary of the domain $\cD$. Treating $\partial\cD$ as a manifold and expressing $\ud K_t$ via $\nabla^2\phi(Y_t)$ and $Z_t$, we discover that $Y$ satisfies the definition of a Brownian $\Gamma$-martingale on the manifold $\partial\cD$ given in \cite{Emery-89}.
 
\begin{Proposition}
\label{prop:variete}
Assume the following:
\begin{itemize}
 \item there exists a convex domain $\mathcal{A}$, satisfying $\bar{\mathcal{A}}\cap \bar{\cD} \subset \partial \cD$,
 \item $\1_{\set{y \in \partial \cD\setminus \bar{\mathcal{A}}}} \nabla d(y,\mathcal{A}) \cdot \nabla \phi(y) \geqslant 0$,
 \item $f\equiv0$ and $\xi \in \bar{\mathcal{A}} \cap \partial\cD$ almost surely,
 \item $(Y,Z,K) \in \SP{2} \times \HP{2} \times \mathscr{K}^1$ solve \eqref{eq reflected bsde}.
\end{itemize}
Then, $Y \in \bar{\mathcal{A}} \cap \bar{\cD} \subset \partial \cD$ almost surely. Moreover, we have
\begin{equation}
\label{eq:vartK}
 \ud \mathrm{Var}_t(K) = \left[-\frac12\mathrm{Tr}[Z_t^\top\nabla^2\phi(Y_t)Z_t] \right]^+ \ud t.
\end{equation}
Finally, $Y$ is a $\Gamma$-martingale, with the prescripted terminal value $\xi$, in the manifold $\partial \cD$ endowed with the Riemannian structured inherited from $\mathbb{R}^d$ and its canonical connection $\Gamma$, as defined in \cite{Emery-89}.
\end{Proposition}

\begin{Remark}
It is worth mentioning that the assumptions made in Proposition \ref{prop:variete} imply that the set $\mathcal{A}$ cannot be smooth. To obtain an intuitive understanding of what the set $\mathcal{A}$ may look like, we refer the reader to the example that follows.
\end{Remark}

\proof
We apply It\^o's formula for general convex functions (in the form of an inequality, as in \cite{Bouleau-84}) to the process $d(Y_t,\mathcal{A})$ to obtain
\begin{align*}
 0 \leqslant d(Y_t,\mathcal{A}) \leqslant \mathbb{E}_t\left[d(\xi,\mathcal{A})-\int_t^T \1_{\set{Y_s \in \partial \cD\setminus \bar{\mathcal{A}}}} \nabla d(Y_s,\mathcal{A}) dK_s\right] \leqslant 0, \quad t \in [0,T],
\end{align*}
which gives us $Y \in \bar{\mathcal{A}} \cap \cD \subset \partial \cD$. Applying It\^o's formula to $\phi(Y_t)$ yields \eqref{eq:vartK}.
Finally, using \eqref{eq:vartK}, the fact that $dK_t$ is orthogonal to the tangent space of $\partial \cD$ at the point $Y_t$, as well as (4.9), (4.10), (5.6)(ii) from \cite{Emery-89}, we conclude that $Y$ is a $\Gamma$-martingale on $\partial\cD$. 
\eproof

\smallskip

In the remainder of this section, we assume $f=0$ and present a simple example of the domains $\cD$ and $\mathcal{A}$ for which the assumptions of Proposition \ref{prop:variete} hold. This example allows us to obtain an alternative proof of a known result on $\Gamma$-martingales using the reflected BSDEs and, on the other hand, to illustrate the sharpness of the weak star-shape assumption (see Assumption \ref{ass:intro.domain}) using the martingales on manifolds.

In this example, we construct the functions $\phi$ and $\phi_{\cC}$, which define the domains $\cD$ and $\mathcal{C}$ as in Assumption \ref{ass:intro.domain}, first, in the plane $\mathcal{P}:=\R\times\{0\}^{d-2}\times \R$ of $\R^d$. {These functions are built symmetric with respect to the $y_d$-axis, as the domains are themselves symmetric, see the precise description below}. Then, we extend them to $\R^d$ via
$$
\phi(y) = \phi((r(y),0,...,0,y_d)),\quad \phi_{\mathcal{C}}(y) = \phi_{\mathcal{C}}((r(y),0,...,0,y_d)),
$$
with $r(y):=\left(\sum_{i=1}^{d-1} |y^i|^2\right)^{1/2}$. By an abuse of notation we denote by the same names the associated domains $\cD$ and $\mathcal{C}$ constructed in $\mathbb{R}^d$ and their intersections with $\mathcal{P}$.

{
We consider the three parameters $\alpha \in (0,\pi/2)$, $\eta>0$, $\varepsilon\in(0,\pi/2-\alpha)$, and the domains $\cD_{\alpha,\eta,\varepsilon}$, $\cC_{\alpha,\eta}$, $\mathcal{A}_{\alpha,\varepsilon}$ given in Figure \ref{figure}, which satisfy the following properties. 
\begin{itemize}
 \item $\mathcal{C}_{\alpha,\eta}$ is a square centered at $(0,-1-\eta-\sin(\alpha))$, with length $2\sin(\alpha)+2\eta$, with the edges parallel to axes and with rounded corners (obtained by modifying the square in the $\eta$-neighborhoods of these corners), such that $\partial \mathcal{C}_{\alpha,\eta}$ is a $C^2$ curve and $\mathcal{C}_{\alpha,\eta}$ is convex, 
 \item $\cD_{\alpha,\eta,\varepsilon}$ is symmetric w.r.t. the axis $y_d$.
 \item $\partial \cD_{\alpha,\eta,\varepsilon}$ is $C^2$ and is made up of the following pieces: 
 \begin{itemize}
 \item the arc $\mathcal{S}_{\alpha}$ of angle $2\alpha$, symmetric w.r.t. the axis $y_d$, of the circle centered at zero and with the radius $1$, 
 \item the arc of angle $2\alpha$, symmetric w.r.t. the axis $y_d$, of the circle centered at zero and with the radius $(2\sin(\alpha)+2\eta+1)/\cos(\alpha)$,  
 \item and two smooth curves $\mathcal{L}_{1}$ and $\mathcal{L}_{2,}$, symmetric to each other w.r.t. the axis $y_d$, which connect the two arcs described above forming a $C^2$ closed curve that does not intersect itself nor $\cC_{\alpha,\eta}$.
 \end{itemize}
 \item We denote by $A^1$ (respectively, $A^2$) the end point of the curve $\mathcal{S}_{\alpha}$ that belongs to the right (respectively, left) half-plane w.r.t. the axis $y_d$.
 \item Let us assume that $\mathcal{L}_{1}$ (respectively, $\mathcal{L}_{2}$) belongs to the right (respectively, left) half-plane w.r.t. the axis $y_d$. We also assume that the curve $\mathcal{L}_{1}$ is constructed so that, in its natural parameterization with the starting point $A^1$, the slope of its tangent vector has exactly one change of monotonicity. Namely, we assume that there exists a point $B^1_{\varepsilon}$, such that the angle between $B^1_{\varepsilon}$ and $A^1$ relative to zero is $\varepsilon$ and such that the derivative of the slope of the aforementioned tangent vector is continuous, nonincreasing, and equal to zero at $B^1_{\varepsilon}$.
The curve $\mathcal{L}_{2}$, then, satisfies the analogous property due to symmetry, with the associated point $B^2_{\varepsilon}$.
\item As the curve $\partial \cD_{\alpha,\eta,\varepsilon}$ is $C^2$, closed, and without self-intersections, we construct $\phi$ as the signed distance to $\partial \cD_{\alpha,\eta,\varepsilon}$ in a neighborhood of $\partial \cD_{\alpha,\eta,\varepsilon}$ and, then, extend it in a smooth way to $\mathbb{R}^2$. $\phi_{\cC}$ is constructed similarly.
 \item We define $\mathcal{S}_{\alpha,\varepsilon}$ as the concatenation of the curves $\overset{\frown}{B^2_{\varepsilon}A^2}$, $\mathcal{S}_{\alpha}$, $\overset{\frown}{A^1B^1_{\varepsilon}}$, and we define $\mathcal{A}_{\alpha,\varepsilon}$ as the interior of the convex hull of $\mathcal{S}_{\alpha,\varepsilon}$.
 \item Finally, we always assume that $\eta>0$ is small enough, such that $\mathcal{C}_{\alpha,\eta}$ is included in the triangle with vertices $P^1$, $0$ and $P^2$, as shown in Figure \ref{figure}. This ensures that $\bar{\mathcal{C}}_{\alpha,\eta}\subset \cD_{\alpha,\eta,\varepsilon}$ for any $\varepsilon>0$.
 \end{itemize}

\usetikzlibrary{math} 

\tikzmath{ \et = 0.4; \eps = 15; \alph = 45; \lamb = sin(\alph)+\et; \t={(sin(2*\alph)-sin(2*(\alph+\eps)))/(cos(2*\alph)-cos(2*(\alph+\eps)))}; \tp={(sin(2*\alph)-sin(2*(\alph+\eps/2)))/(cos(2*\alph)-cos(2*(\alph+\eps/2)))}; \coeff=5;}

 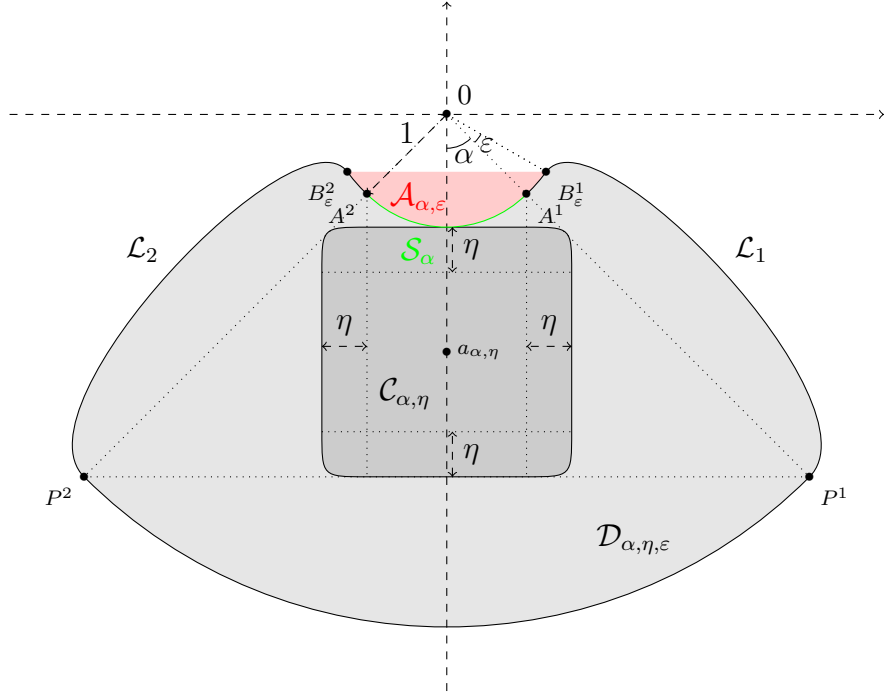
\begin{figure}
 \begin{center}
\begin{tikzpicture}[scale = 1.5]
 
 \coordinate (Op) at (0,\lamb+1);
 \coordinate (O) at (0,0);
 \coordinate (A1) at ({sin(\alph)},{\lamb+1 -cos(\alph)});
 \coordinate (A2) at ({-sin(\alph)},{\lamb+1 -cos(\alph)});
 \coordinate (B1) at 
 ( {0 + (sin(\alph)+\t*cos(\alph)+sin(\alph+\eps))/2} , {\lamb+1 +(-cos(\alph)+\t*sin(\alph)-cos(\alph+\eps))/2} ) ;
 \coordinate (B1aux) at 
 ( {0 + sin(\alph)+\tp*cos(\alph)} , {\lamb+1 -cos(\alph)+\tp*sin(\alph)} ) ;
 \coordinate (B1ter) at ({(1+\coeff*\lamb)*(sin(\alph)+\t*cos(\alph)+sin(\alph+\eps))/2-\coeff*\lamb*(sin(\alph)+\tp*cos(\alph))},{(1+\coeff*\lamb)*(\lamb+1 +(-cos(\alph)+\t*sin(\alph)-cos(\alph+\eps))/2)- \coeff*\lamb*(\lamb+1 -cos(\alph)+\tp*sin(\alph))}) ;
 \coordinate (B2) at 
 ( {0 - (sin(\alph)+\t*cos(\alph)+sin(\alph+\eps))/2} , {\lamb+1 +(-cos(\alph)+\t*sin(\alph)-cos(\alph+\eps))/2} ) ;
 \coordinate (B2aux) at 
 ( {0 - sin(\alph)-\tp*cos(\alph)} , {\lamb+1 -cos(\alph)+\tp*sin(\alph)} ) ;
 \coordinate (B2ter) at ({-(1+\coeff*\lamb)*(sin(\alph)+\t*cos(\alph)+sin(\alph+\eps))/2+\coeff*\lamb*(sin(\alph)+\tp*cos(\alph))},{(1+\coeff*\lamb)*(\lamb+1 +(-cos(\alph)+\t*sin(\alph)-cos(\alph+\eps))/2)- \coeff*\lamb*(\lamb+1 -cos(\alph)+\tp*sin(\alph))}) ;
 \coordinate (C1) at ({sin(\alph)*(1+2*\lamb)/cos(\alph)},{\lamb+1 -(1+2*\lamb)*cos(\alph)/cos(\alph)});
 \coordinate (C1aux) at (({(sin(\alph)+0.2*cos(\alph))*(1+2*\lamb)/cos(\alph},{\lamb+1 +(1+2*\lamb)*(-cos(\alph)+0.2*sin(\alph))/cos(\alph)});
 \coordinate (C2) at ({-sin(\alph)*(1+2*\lamb)/cos(\alph)},{\lamb+1 -(1+2*\lamb)*cos(\alph)/cos(\alph)});
 \coordinate (C2aux) at (({-(sin(\alph)+0.2*cos(\alph))*(1+2*\lamb)/cos(\alph)},{\lamb+1 +(1+2*\lamb)*(-cos(\alph)+0.2*sin(\alph))/cos(\alph)});

 \fill[color=black!10] (A1) .. controls (B1aux) and (B1aux) .. (B1) .. controls (B1ter) and (C1aux) .. (C1) arc  (-90+\alph:-90-\alph:{(1+2*\lamb)/cos(\alph)}) .. controls (C2aux) and (B2ter) .. (B2) .. controls
 (B2aux) and (B2aux) .. (A2) arc ((-90-\alph:-90+\alph:1) -- cycle ;
 \draw (1.5*\lamb,-1.5*\lamb) node {$\mathcal{D}_{\alpha,\eta,\varepsilon}$} ;

 \fill[color=black!20] (-\lamb+\et,\lamb) -- (\lamb-\et,\lamb) .. controls (\lamb,\lamb) and (\lamb,\lamb) .. (\lamb,\lamb-\et) -- (\lamb,-\lamb+\et) .. controls (\lamb,-\lamb) and (\lamb,-\lamb) .. (\lamb-\et,-\lamb) -- (-\lamb+\et,-\lamb) .. controls (-\lamb,-\lamb) and (-\lamb,-\lamb) .. (-\lamb,-\lamb+\et) -- (-\lamb,\lamb-\et) .. controls (-\lamb,\lamb) and (-\lamb,\lamb) .. (-\lamb+\et,\lamb);
 \draw (-\lamb+\et,\lamb) -- (\lamb-\et,\lamb) .. controls (\lamb,\lamb) and (\lamb,\lamb) .. (\lamb,\lamb-\et) -- (\lamb,-\lamb+\et) .. controls (\lamb,-\lamb) and (\lamb,-\lamb) .. (\lamb-\et,-\lamb) -- (-\lamb+\et,-\lamb) .. controls (-\lamb,-\lamb) and (-\lamb,-\lamb) .. (-\lamb,-\lamb+\et) -- (-\lamb,\lamb-\et) .. controls (-\lamb,\lamb) and (-\lamb,\lamb) .. (-\lamb+\et,\lamb);
 \draw (-\lamb/3,-\lamb/3) node {$\mathcal{C}_{\alpha,\eta}$} ;
 
 \fill[color=red!20] (A1) .. controls (B1aux) and (B1aux) .. (B1) -- (B2) .. controls
 (B2aux) and (B2aux) .. (A2) arc ((-90-\alph:-90+\alph:1) -- cycle ;
 \draw[red] (-0.25,\lamb+0.25) node {$\mathcal{A}_{\alpha,\varepsilon}$} ;

 \draw[->,dashed] (0,-\lamb-\et-1.5) -- (0,\lamb+2) ;
 \draw[->,dashed] (-3.5*\lamb,\lamb+1) -- (3.5*\lamb,\lamb+1) ;

 \draw[green] (A2) arc (-90-\alph:-90+\alph:1)node[midway,below left] {$\mathcal{S}_{\alpha}$} ;

 \draw (A1) .. controls (B1aux) and (B1aux) .. (B1) ;
 \draw (B1) .. controls (B1ter) and (C1aux) .. (C1) node [midway,above right] {\small $\mathcal{L}_1$};
 \draw (A2) .. controls (B2aux) and (B2aux) .. (B2) ;
 \draw (B2) .. controls (B2ter) and (C2aux) .. (C2) node [midway,above left] {\small $\mathcal{L}_2$};
 
 \draw (A1) node {\tiny $\bullet$};
 \draw (A1) node [below right] {\tiny $A^1$};
 \draw (A2) node {\tiny $\bullet$};
 \draw (A2) node [below left] {\tiny $A^2$};
 \draw (B1) node {\tiny $\bullet$};
 \draw (B1) node [below right] {\tiny $B^1_{\varepsilon}$};
 \draw (B2) node {\tiny $\bullet$};
 \draw (B2) node [below left] {\tiny $B^2_{\varepsilon}$};
 \draw (C1) node {\tiny $\bullet$};
 \draw (C1) node [below right] {\tiny $P^1$};
 \draw (C2) node {\tiny $\bullet$};
 \draw (C2) node [below left] {\tiny $P^2$};
 
 \draw (C1) arc  (-90+\alph:-90-\alph:{(1+2*\lamb)/cos(\alph)});
 
 \draw (O) node {\tiny $\bullet$};
 \draw (O) node [right] {\tiny $a_{\alpha,\eta}$};
 \draw[dotted] (Op) -- (C1) ;
 \draw[dotted] (Op) -- (C2) ;
 \draw[dashed,<->] (Op) -- (A2) node[midway,above] {$1$};
 \draw[dashed,<->] (0.05,\lamb) -- (0.05,\lamb-\et) node[midway,right] {$\eta$};
 \draw[dashed,<->] (\lamb,0.05) -- (\lamb-\et,0.05) node[midway,above] {$\eta$};
 \draw[dashed,<->] (0.05,-\lamb) -- (0.05,-\lamb+\et) node[midway,right] {$\eta$};
 \draw[dashed,<->] (-\lamb,0.05) -- (-\lamb+\et,0.05) node[midway,above] {$\eta$};
 \draw (0,\lamb+0.7) arc (-90:-90+\alph:0.3) ;
 \draw ({0.4*sin(\alph/2)},{\lamb+1-0.4*cos(\alph/2)}) node {$\alpha$} ;
 \draw[dotted] (Op) -- (B1) ;
 \draw (({0.35*cos(\alph)},{\lamb+1-0.35*sin(\alph)}) arc (-90+\alph:-90+\alph+\eps:0.35) ;
 \draw ({0.45*sin(\alph+\eps/2)},{\lamb+1-0.45*cos(\alph+\eps/2)}) node {$\varepsilon$} ;
 \draw[dotted] (Op) -- (B1) ;
 \draw (0,1+\lamb) node {\tiny $\bullet$};
 \draw (0,1+\lamb) node [above right] {\small$0$};
 \draw[dotted] (A1) -- (\lamb-\et,-\lamb);
 \draw[dotted] (A2) -- (-\lamb+\et,-\lamb);
 \draw[dotted] (-\lamb,\lamb-\et) --(\lamb,\lamb-\et);
 \draw[dotted] (-\lamb,-\lamb+\et) --(\lamb,-\lamb+\et);
 \draw[dotted] (C2) -- (C1);
 
\end{tikzpicture}
 \end{center}
 \caption{Domains $\cD_{\alpha,\eta,\varepsilon}$, $\cC_{\alpha,\eta}$, $\mathcal{A}_{\alpha,\varepsilon}$}
  \label{figure}
 \end{figure}

Let us now consider a terminal condition $\xi \in \mathcal{S}_{\alpha}$ and verify that $\cD_{\alpha,\eta,\varepsilon}$, $\mathcal{A}_{\alpha,\varepsilon}$ and $\xi$ satisfy the desired assumptions. We easily deduce that $R_0=1$. Then, for any $\alpha \in (0,\pi/2)$ and $\eta>0$, there exists $\varepsilon_0\in(0,\pi/2-\alpha)$, such that, for all $0<\varepsilon<\varepsilon_0$, the condition \eqref{def:gamma} holds up to a shift of coordinates in $\R^d$ that maps the origin to $a_{\alpha,\eta}:=(0,...,0,-1-\eta-\sin(\alpha))$. The other conditions of Assumption \ref{ass:intro.domain} follow easily. 

Next, we notice that, in the discrete path-dependent framework and under Assumption \ref{ass:pathdependent}, we can apply Theorem \ref{thm:pathdependent} to conclude that there 
exists a unique (in $\mathscr{U}(1)$) triplet $(Y^{\varepsilon},Z^{\varepsilon},K^{\varepsilon})\in\SP{2}\times \HP{2} \times \mathscr{K}^1$ that solves \eqref{eq reflected bsde} in the domain $\cD=\cD_{\alpha,\eta,\varepsilon}$ (we suppress the dependence of the solution on $\eta$ and $\alpha$ as  they are fixed in what follows).
Proposition \ref{prop:variete} applied to $\cD=\cD_{\alpha,\eta,\varepsilon}$ and $\mathcal{A}=\mathcal{A}_{\alpha,\varepsilon}$ (whose assumptions are satisfied by the construction of $\mathcal{A}_{\alpha,\varepsilon}$, $\mathcal{L}_1$ and $\mathcal{L}_2$) states that $Y^{\varepsilon}$ lives in $\mathcal{S}_{\alpha,\varepsilon}$. Then, the stability result of Proposition \ref{prop:stability2} yields that $\{(Y^{1/n},Z^{1/n},K^{1/n})\}_{n=1}^{\infty}$ is a Cauchy sequence and, hence, has a limit $(Y,Z,K)$. It is clear that $Y$ lives on $\mathcal{S}_{\alpha}$. Then, applying the arguments similar to those used in the proof of Theorem \ref{thm:existence:Markov}, one can deduce that $(Y^{},Z^{},K^{})$ solve the reflected BSDE \eqref{eq reflected bsde} in the domain $\cD=\cD_{\alpha,\eta,\varepsilon'}$, for any $\varepsilon'\in(0,\varepsilon_0)$. Applying Proposition \ref{prop:variete} once more and recalling that $Y^{}$ lives on $\mathcal{S}_{\alpha}$, we conclude that $Y^{}$ is a $\Gamma$-martingale on the manifold $\mathcal{S}_{\alpha}$ with the terminal condition $\xi$. The uniqueness part of Theorem \ref{thm:pathdependent} yields that such a $\Gamma$-martingale is unique (in $\mathscr{U}(1)$). 

Now, let us come back to a general terminal condition. We first notice that Proposition \ref{A priori estimate reflected BSDE} holds for any solution $(Y,Z,K) \in \SP{2}\times \HP{2}\times \mathscr{K}^1$ of \eqref{eq reflected bsde} that lives in $\mathcal{S}_{\alpha}$ and satisfies Assumption \ref{ass:general}(i) with $\gamma$ replaced by 
$$\gamma_{\alpha} : = \inf_{y \in \mathcal{S}_{\alpha}} \nabla \phi_{\mathcal{C}}(y)  \cdot \frac{\nabla \phi(y)}{|\nabla \phi(y)|}.$$
We can easily compute $\gamma_{\alpha} = \cos(\alpha)$. Moreover, we also have $R_0=1$ and 
$$|\phi_{\mathcal{C}}^+(\xi)|_{\mathscr{L}^{\infty}} \leqslant 1-\cos(\alpha).$$
Thus, we conclude that Assumption \ref{ass:general}(i) is fulfilled with $\theta=2$ as long as $\cos(\alpha)>2/3$. 
Considering a sequence of discrete path-dependent terminal conditions that approximate the given (general) terminal condition and take values in $\mathcal{S}_{\alpha}$, we repeat the proof of Theorem \ref{thm:general} obtaining a unique (in $\mathscr{U}(2)$) triplet $(Y^{},Z^{},K^{})\in\SP{2}\times \HP{2} \times \mathscr{K}^1$ that solves \eqref{eq reflected bsde} in the domain $\cD=\cD_{\alpha,\eta,\varepsilon'}$, for any $\varepsilon'\in(0,\varepsilon_0)$, and is such that $Y^{}$ lives in $\mathcal{S}_{\alpha}$. Applying Proposition \ref{prop:variete} once more, we conclude that $Y^{}$ is a $\Gamma$-martingale on the manifold $\mathcal{S}_{\alpha}$ with the terminal condition $\xi$. The uniqueness part of Theorem \ref{thm:pathdependent} yields that such a $\Gamma$-martingale is unique in $\mathscr{U}(2)$.

\smallskip

To sum up, the above construction proves the existence and uniqueness result for a Brownian $\Gamma$-martingale with a prescribed discrete path-dependent terminal condition $\xi$, satisfying Assumption \ref{ass:pathdependent}, on any sector of the sphere $\mathbb{S}^{d-1}$ (we understand a sector as an intersection of a sphere and a half-space) that is strictly contained in a hemisphere. For a general terminal condition $\xi$, we are only able to tackle the case $\alpha<\text{arccos}(2/3)$. Thus, our results provide an alternative proof of a particular case of \cite{Kendall-90,Picard-91}, where the existence and uniqueness for any $\alpha<\pi/2$ is established. Considering the case $\alpha=\pi/2$, we notice that, for any $\cD$ that is included in the outside of a sphere and whose boundary contains a hemisphere, it is impossible to find a convex domain $\cC \subset \cD$ that can ``see'' all points on the boundary of this hemisphere with a strictly positive angle: in other words, \eqref{def:gamma} can not be fulfilled. 
Therefore, our existence and uniqueness results fail for the case of a hemisphere. On the other hand, considering directly the problem of existence and uniqueness of a Brownian $\Gamma$-martingale with a prescribed terminal condition on a closed hemisphere of $\mathbb{S}^{d-1}$, we notice the major challenge that is due to the non-uniqueness of geodesics, when $d\geqslant 3$. Indeed, assume that $\xi$ takes its values in the set $\{z_1,z_2\}$ consisting of two antipodes on the sphere (i.e., the line connecting the two points goes through the center of the sphere) and note that $\mathcal{S}_{\pi/2}$ does contain such points. Then, for any shortest arc $\overset{\frown}{z_1z_2} \subset \mathcal{S}_{\pi/2}$, there exists  a $\Gamma$-martingale on the manifold $\overset{\frown}{z_1z_2}$ with the terminal condition $\xi$. As the arc $\overset{\frown}{z_1z_2}$ is a geodesic, we conclude that the resulting $\Gamma$-martingale is also a $\Gamma$-martingale in the larger manifold $\mathcal{S}_{\pi/2}$. Assuming that $\xi$ takes each of its two values with strictly positive probability and recalling that there are infinitely many geodesic arcs $\overset{\frown}{z_1z_2}$ on $\mathcal{S}_{\pi/2}$, we conclude that the uniqueness of a $\Gamma$-martingale on $\mathbb{S}^{d-1}$ with the terminal condition $\xi$ does not hold. This observation, in particular, illustrates the sharpness of  the weak star-shape property in Assumption \ref{ass:intro.domain} (condition \eqref{def:gamma}) in the case of a general terminal condition and general $d\geq2$.

\smallskip

Let us also mention that the non-uniqueness described above does not occur for $d=2$, which indicates that it may be possible to relax our assumptions for reflected BSDEs in planar non-convex domains. In particular, we refer to \cite{Picard-89} for a complete treatment of $\Gamma$-martingales on $\mathbb{S}^{1}$. The latter result also yields the existence and uniqueness of a solution to the reflected BSDE in the domain $\cD =\{y \in \mathbb{R}^2, 1<|y|<2\}$, which does not possess the weak star-shape property, with zero generator and with a terminal condition satisfying $|\xi|=1$.

Moreover, in Section $3$ of \cite{Picard-91}, Picard was able to prove the existence and uniqueness of a Brownian $\Gamma$-martingale with a prescribed terminal condition in a closed hemisphere of $\mathbb{S}^{d-1}$, and in an even bigger domain, for a small enough $T$ and under a smoothness assumption on the terminal condition\footnote{To be precise, it is assumed that the process $Z$, defined by $\xi=\mathbb{E}[\xi]+\int_0^T Z_s \ud W_s$, has sufficiently small $\int_0^T \text{ess sup}_{\Omega} |Z_s|^2\ud s$.}.
The latter indicates that in a smooth Markovian or discrete path-dependent framework, under an additional smallness assumption, it may also be possible to relax the requirement of a weak star-shape property even for $d>2$.

\smallskip

Finally, let us give a simple example showing that the \emph{a priori} estimates of Proposition \ref{A priori estimate reflected BSDE} are not sharp.\footnote{Note that these estimates are not needed in a Markovian or discrete path-dependent case.} Mimicking \cite{Picard-89}, we consider a $\mathcal{F}_T$-measurable random variable $\nu$ with values in $[-\alpha,\alpha]$, where $0<\alpha<\pi/2$ is a given parameter, and let $(\theta_t,\eta_t)_{t \in [0,T]}$ be the solution of the BSDE $\theta_t = \nu -\int_t^T \eta_s \ud W_s$ for $t\in [0,T]$. We set $Y_t = (\cos (\theta_t),\sin (\theta_t))^\top$ for all $t \in [0,T]$ and we easily check that $Y$ is a solution to the BSDE
$$Y_t = \xi +\int_t^T \frac{|Z_s|^2}{2}Y_s \ud s -\int_t^T Z_s \ud W_s, \quad 0 \leqslant t \leqslant T,$$
where $\xi = (\cos (\nu),\sin (\nu))^\top$ and $Z_t = (-\eta_t\sin (\theta_t),\eta_t \cos (\theta_t))^\top$. Notice that this multidimensional quadratic BSDE can also be seen as a reflected BSDE in the domain
$\cD_{\alpha,\eta,\varepsilon}$, with sufficiently small $\eta,\varepsilon>0$, rotated by $\pi/2$. Indeed, $Y$ lives on (rotated) $\mathcal{S}_\alpha$ and the reflecting term is always pointing along the outer normal vector to (rotated) $\mathcal{S}_\alpha$. Recall that $\cD_{\alpha,\eta,\varepsilon}$ satisfies the weak star-shape property and note that $\ud \text{Var}_t(K) = \frac{|\eta_t|^2}{2}$. Then, an application of It\^o's formula to $\theta_t^2$ yields
$$
\mathbb{E}_t\left[\int_t^T \ud \text{Var}_s(K)\right] = \frac{1}{2} \mathbb{E}_t\left[ \int_t^T |\eta_s|^2\ud s\right] =
\frac{1}{2} \mathbb{E}_t\left[ \nu^2 - (\EE_t\nu)^2\right] \leqslant \frac{\alpha^2}{2}.
$$
Moreover, the above becomes an equality for $t=0$ and $\nu =  \text{sign}(W_T)\alpha$.
Then, recalling that $R_0=1$ for $\cD_{\alpha,\eta,\varepsilon}$, we deduce from the John-Nirenberg inequality that 
$$ \mathbb{E}\left[ e^{\frac{2 p}{R_0} \mathrm{Var}_T(K)} \right] < \infty,$$
for some $p>1$, provided $\alpha<1$, which is weaker than the condition $\alpha<\text{arccos}(2/3)<1$ required by Assumption \ref{ass:general}(i) with $\theta=2$, as computed earlier in this subsection. 
}

\bibliographystyle{plain}
\def\cprime{$'$}

\end{document}